\pdfoutput=1
\documentclass[11pt,a4paper]{article}

\usepackage[T1]{fontenc}
\usepackage[utf8]{inputenc}
\usepackage{amsmath,amssymb,amsthm}
\usepackage{graphicx}
\usepackage{booktabs}
\usepackage{caption}
\usepackage[margin=2.6cm]{geometry}

\newcommand{\TN}{the technical note~\cite{repo}}

\usepackage{xcolor}
\usepackage{hyperref}
\usepackage{xcolor}
\usepackage{textcomp}   %

\hypersetup{colorlinks=true,linkcolor=blue!50!black,citecolor=blue!50!black,
            urlcolor=blue!50!black}

\theoremstyle{definition}

\newcommand{\wmax}{V}

\title{\bf GPU-accelerated wind farm layout search\\
       with a distilled endogenous wake model (EndoWake)}
\author{%
Martina Fischetti\thanks{Department of Statistics and Operations Research
and Institute of Mathematics (IMUS), University of Seville, Facultad de
Matem\'aticas, Av.~Reina Mercedes s/n, 41012 Sevilla, Spain.
ORCID: 0000-0002-7673-6917.}
\and
Matteo Fischetti\thanks{Department of Information Engineering, University of
Padova, Via Gradenigo 6/B, 35131 Padova, Italy.
ORCID: 0000-0001-6601-0568.
\emph{Corresponding author}: \texttt{matteo.fischetti@unipd.it}}}

\date{}

\begin{document}
\maketitle
\begin{abstract}
We build on EndoWake, an endogenous wake model that links the
wind speed at every grid cell to its upwind neighbors by linear
constraints, so that wake field, siting decisions and project constraints
share one mixed-integer model. Its four parameters are calibrated once on
the single wake of a Reynolds-averaged Navier--Stokes (RANS) simulation,
and a graphics processing unit (GPU) scores over a million layouts per
second per wind direction. At Lillgrund, we compare searches guided by
EndoWake and by PyWake under a RANS judge that no search calls: the
fidelity of a wake model is not the quality of the layouts optimized
under it. Two diagnostics separate a directional quadrature error between
the sampled wind directions from a spatially localized surrogate bias off
each wake edge, not
reported before to our knowledge; at random directions, five paired
starts detect no difference between the pipelines. The main result
puts the speed to use: EndoWake ranks the moves of a local search on the
GPU, and PyWake confirms every accepted move. On twenty unseen starts,
with one GPU and thirty CPU processes, screen-and-confirm reaches
similar mean gains to a search on PyWake alone,
7.9\ times
faster at one wind speed and 14.9\ times
faster on the annual energy (3.0\
and 4.7%
\ with a proof of local optimality
charged to the timed search); on a second site it keeps its speed,
and its local optima fall slightly short of those of PyWake alone. Code
and data are deposited on Zenodo.
\end{abstract}


\section{Introduction}
\label{sec:intro}

Wind farm layout optimization places the turbines of a farm within an
awarded area so as to maximize its annual energy production, under the
conditions that come with the award: an installed capacity, hence a maximum
number of machines, no-go areas, position-dependent costs, and obstacles
that shelter the flow~\cite{fischetti2019overview,fischetti2020vattenfall}.
Over a concession of twenty-five years or so, a fraction of a percent of
annual energy is a sum worth contesting, and the bottleneck of the search
is the wake computation. Engineering models, such as Jensen's
cone~\cite{jensen1983} or the Gaussian wake of Bastankhah and
Port\'e-Agel~\cite{bastankhah2014}, are approximate but fast, hence they are
the models that optimization frameworks call millions of times.
Computational fluid dynamics (CFD) is the reference against which they are
validated~\cite{pena2018anholt,sebastiani2021} but, adjoint methods
aside~\cite{king2017,antonini2020adjoint}, too slow to sit inside a search:
the Reynolds-averaged Navier--Stokes (RANS) reference of this paper costs
405~s per wind direction for one layout.

Engineering models drive three kinds of layout optimizer: a discrete site
selection on an interference
matrix~\cite{turner2014,zhang2014,fischetti2016proximity,%
fischetti2021nextgen,fischetti2023integrated,perezrua2023}, a
metaheuristic~\cite{cazzaro2022vns,cazzaro2022multiscale}, or a gradient on
continuous coordinates, as in TOPFARM~\cite{topfarm} on top of
PyWake~\cite{pywake}, the open-source simulator of DTU Wind Energy, in
FLORIS~\cite{floris}, and in FLOWERS~\cite{locascio2022,locascio2024,%
locascio2025}. Position-dependent costs and exclusion
zones~\cite{criadorisco2024zones} can be written on continuous coordinates
too, whereas counting constraints (cardinalities, zone caps) call for a
discrete treatment. On binary site-selection variables, all these
conditions are linear constraints, which a wake model linear in the same
variables can join in one formulation.

\begin{figure}[t]
\centering
\includegraphics[width=\textwidth]{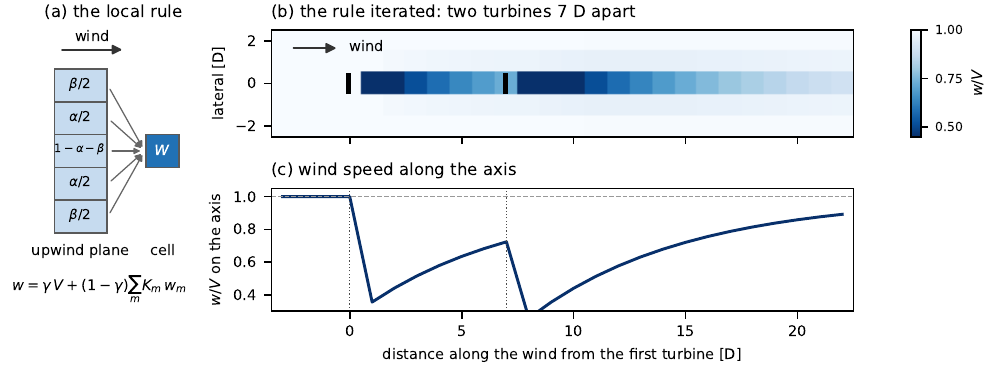}
\caption{The endogenous field, with the kernel distilled from our RANS at
a cell of one rotor diameter $D$. (a)~The local rule: the speed at a cell
is $\gamma$ times the free stream plus $1-\gamma$ times a weighted average
of five upwind cells, with weights $(\beta/2, \alpha/2, 1-\alpha-\beta,
\alpha/2, \beta/2)$; behind a turbine, the speed is capped at a fixed
fraction of the speed the turbine sees. (b)~The rule iterated for two
turbines $7D$ apart along the wind. (c)~Speed along the axis: the width of
the wake, its recovery and the deeper deficit behind the second turbine
emerge from the repetition.}
\label{fig:intro-rule}
\end{figure}

\paragraph{EndoWake: wind propagation as an endogenous field.}
EndoWake, the wake model that we introduced in~\cite{companion}, is such a
model. It treats the wind speed at every cell of a grid as a continuous
decision variable, linked to its upwind neighbors by linear propagation
constraints, and a turbine multiplicatively reduces the speed it sees
(Figure~\ref{fig:intro-rule}). Wake superposition, higher-order interactions
and wake recovery all follow from iterating that local rule through the
grid, and we call the result the \emph{endogenous field}. The rule has four
numbers ($\alpha$, $\beta$, $\gamma$ and the source amplitude), calibrated
once against the deficit field of one isolated turbine: this is what a heavy
CFD model is distilled into in Section~\ref{sec:distill}. Other wake models
replace the ad hoc superposition rule by momentum
conservation~\cite{zong2020}, by a merging method~\cite{lanzilao2022} or by
a cumulative formulation~\cite{bastankhah2021,zhang2026}, but they remain
functions that the optimizer calls from outside. In EndoWake, the
propagation is written \emph{inside} the optimization model, as linear
constraints that couple the wind speeds to the binary site-selection
variables of a mixed-integer linear programming (MILP) model.

\paragraph{Fidelity and decisions.}
A search looks explicitly for the spots where its objective is largest,
hence small, localized optimistic errors of an otherwise precise model are
found and exploited---the optimizer's curse~\cite{smith2006}. That the
fidelity of a model and the quality of the layouts optimized under it can
part ways has been observed in the wind literature
too~\cite{thomas2023les,zhang2026}: in a recent comparison of
RANS-calibrated engineering models with a neural surrogate, the model
of lowest prediction error did not give the layouts of highest
RANS-validated energy~\cite{schoeler2026validation}. A discretized layout objective carries
two such errors, one per discretization. The \emph{directional quadrature
error} comes from sampling the wind rose at a finite set of directions:
between two sampled directions a pair of machines looks out of each
other's wake, and a search is drawn to these gaps. The effect is known to
practitioners~\cite{fischetti2017thesis,feng2015,quaeghebeur2021}; here we
measure it direction by direction against a RANS judge. The
\emph{spatially localized surrogate bias} comes from the discretization of
the area in a propagator written on a grid: at a cell of $D/3$ the EndoWake
field is more faithful on average, but it is systematically optimistic in
a band just off the edge of each wake. To our knowledge this bias has not
been reported before; it belongs to the discretized propagator, since an
analytical model evaluated at the same coordinates returns the same field
whatever the candidate sites.

\paragraph{What this paper asks.}
The question is what EndoWake is worth to a planner, and the answer has two
parts. First, a more faithful field does not by itself give better layouts:
refining the EndoWake field brings it closer to the RANS wake but not to
better layouts, a search guided by EndoWake in place of PyWake is drawn to
the errors of the discretization, and at our judge five paired starts
detect no difference between the two. Second, the throughput of
EndoWake pays when the search is built for it: EndoWake is not asked to
replace the model the planner trusts, it is distilled from that model and
used to decide which moves the trusted model examines. Multi-fidelity
preselection~\cite{jasa2022multifidelity} and the gap between fidelity and
decision quality each have precedents; new is their combination into one
search: a linear propagator distilled from a heavy model, evaluated in
bulk on a GPU, a shortlist sized by the measured costs of the two models,
and a confirmation in the final objective.

\paragraph{What is inherited, and what is new.}
From~\cite{companion} we take the model of Section~\ref{sec:mip}, its
maximal-field property (one forward sweep evaluates a layout exactly), the
MILP formulation with its valid inequalities, and the calibration recipe
on an analytical wake. New in this paper are the following.
\begin{itemize}
\item \emph{Distillation of a CFD field, and a GPU evaluator.} The
calibration needs only the deficit field of one isolated machine, from a
RANS simulation paid once or from an analytical wake model. On the single
wake, the distilled kernel is about as faithful to our RANS reference as
the Gaussian wake fitted on the same field (Section~\ref{sec:rans-calib}).
A batched evaluator scores 1\,093\,103\ layouts per
second for one wind direction on a data-center GPU
(Section~\ref{sec:exp-batch}), and every award condition that is linear in
the binaries enters as one more constraint, on grids of up to
40\,000\ cells (Section~\ref{sec:opt-tool}).
\item \emph{Two diagnostics of a discretized objective, and a comparison
between pipelines.} We measure the directional quadrature error by
direction and the spatially localized surrogate bias by cell, with a RANS
as judge, on the Lillgrund offshore farm, whose machines and row-wise
deficits from supervisory control and data acquisition (SCADA) records are
public. A judge that samples the directions a search optimized on cannot
see the integration error between them: an advantage of the EndoWake
pipeline at the sector centers turns into a deficit halfway between them
(Section~\ref{sec:cones}), and we propose the sensitivity of a result to
where the judge samples the rose as a first control on any tool that
optimizes over a discretized rose. At directions drawn at random, five
paired starts detect no difference between the pipelines, with intervals
that permit materially different outcomes
(Section~\ref{sec:sirens-judge}). The spatially localized bias is measured
on the field itself, independently of any search
(Section~\ref{sec:mermaid}).
\item \emph{Surrogate screening with PyWake confirmation.} The same recipe
distills PyWake itself, and the distillate then screens the moves of a
local search whose objective is PyWake: EndoWake ranks samples of a
neighborhood of up to 153\ thousand moves
on the GPU, and PyWake evaluates the head of the ranking and accepts a move
only if it improves the PyWake value. The length of that list follows from
the measured speeds of the two models, with no parameter tuned on results.
On twenty held-out Lillgrund starts, screening reaches mean PyWake-valued
gains similar to those of a search on PyWake alone, within the paired
intervals of Section~\ref{sec:proxy}, in a search time about an order of
magnitude shorter, at one speed and on the annual energy, and still
several times shorter when the proof of local optimality, a final scan of
the whole neighborhood by PyWake, is charged to the timed search. On a
second site, Horns Rev~1, the screening keeps its speedups, but its
local optima fall short of those of PyWake alone by a small and
consistent margin.
\item \emph{The code.} The replication package~\cite{repo} releases
the local searches, the sector quadrature and the fast evaluator, in its
CUDA form and in the vectorized CPU form that reproduces it, so that every
throughput and every wall-clock time of this paper can be rerun.
\end{itemize}

Sections~\ref{sec:background} to~\ref{sec:opt-tool} recall the wake
models and the layout tools, restate EndoWake, distill the kernel from a
RANS computation and validate the search engine with the operator's
constraints and at industrial size; Section~\ref{sec:sirens} measures the
two discretization errors at the RANS judge; Section~\ref{sec:proxy}
presents the screening search on PyWake, our main result;
Section~\ref{sec:concl} concludes.
Appendices~\ref{app:rans-details}--\ref{app:protocols} give the RANS
reference, the computing environment and the protocols. A technical note
in the replication package~\cite{repo}, the technical note in what
follows, carries the proofs and specifications that the paper only
states, additional tables and the record of the RANS runs. 
\section{State of the art: wake models and layout methods}
\label{sec:background}

In this section we fix the notation on wakes and position our work among
the layout methods closest to it, which Table~\ref{tab:sota} compares.
Guiding a search with a cheaper model has precedents, but none of them
writes the wake field as linear constraints on binary site variables. The
reader is referred to Burton et al.~\cite{burton2011} and to Port\'e-Agel
et al.~\cite{porteagel2020} for the physics.

\paragraph{Notation.} A turbine facing an undisturbed speed $\wmax$ leaves
behind it a region of reduced speed. At a point $x$ diameters downstream
and $y$ diameters off the axis, the \emph{deficit} is $d(x,y) = \wmax -
w(x,y)$, which the engineering models write as a profile that widens and
weakens with $x$: a top-hat spreading at a constant rate in Jensen's
model~\cite{jensen1983,lissaman1979,katic1986}, and a Gaussian whose width
grows linearly at a rate $k$ in Bastankhah and
Port\'e-Agel~\cite{bastankhah2014}. The deficits of turbines shading the
same point are combined by a \emph{superposition rule}, usually their
sum~\cite{niayifar2016} or the square root of the sum of their squares,
and this choice alone changes the calibrated wake-expansion rate by
relative differences on the order of $100\%$~\cite{locascio2023}.
Interference-matrix formulations also take the thrust of each machine at
the free-stream speed, while PyWake~\cite{pywake} propagates from the
effective speed. Our model needs neither approximation, save that a linear
rule is implied where two wakes merge sideways (Section~\ref{sec:mip}).

\subsection{Layout methods}
\label{sec:related}

\begin{table}[tbp]
\centering\footnotesize
\def\rr{\raggedright}
\setlength{\tabcolsep}{3pt}
\caption{Layout methods that approximate or accelerate an expensive
evaluation: the model approximated (and by what), the decision variables,
what is kept exact, and the cost reduced.}
\label{tab:sota}
\begin{tabular}{@{}p{0.15\textwidth}p{0.23\textwidth}p{0.14\textwidth}p{0.24\textwidth}p{0.16\textwidth}@{}}
\toprule
\rr Method & \rr Model approximated & \rr Variables & \rr Kept exact & \rr Cost reduced \tabularnewline
\midrule
\rr Multifidelity, Jasa et al.~\cite{jasa2022multifidelity}
 & \rr a higher-fidelity wake model, by a lower-fidelity one plus a learned correction
 & \rr continuous coordinates
 & \rr higher-fidelity values at the sampled points
 & \rr higher-fidelity calls \tabularnewline
\rr WakeNet~\cite{anagnostopoulos2023wakenet}
 & \rr Gaussian and curl wake models of FLORIS, by a neural network
 & \rr yaw angles, or continuous coordinates
 & \rr nothing exactly; accuracy measured against FLORIS
 & \rr wake evaluations \tabularnewline
\rr FLOWERS~\cite{locascio2022,locascio2024,locascio2025}
 & \rr the rose average of an engineering model, by a closed form
 & \rr continuous coordinates
 & \rr gradients of its own model
 & \rr integration over the rose \tabularnewline
\rr Fast gradient search~\cite{padron2019pce,quick2023sgd,valottarodrigues2024speeding}
 & \rr none; the rose is integrated by quadrature or sampled
 & \rr continuous coordinates
 & \rr the engineering model
 & \rr evaluations per iteration; wall-clock time \tabularnewline
\rr Interference matrix~\cite{turner2014,fischetti2016proximity,cazzaro2022vns}
 & \rr an engineering model, by precomputed pairwise deficits
 & \rr binary site selection
 & \rr counting and logical constraints; MILP bounds
 & \rr wake evaluation during the search \tabularnewline
\rr EndoWake, Sections~\ref{sec:surrogate}--\ref{sec:sirens}
 & \rr an isolated wake (RANS or an engineering model), by a linear propagator on the grid
 & \rr binary site selection
 & \rr linear side constraints; MILP certificate on small grids
 & \rr CFD, called once offline \tabularnewline
\rr Screen-and-confirm, Section~\ref{sec:proxy}
 & \rr PyWake, by EndoWake fitted to its single wake (ranking only)
 & \rr binary site selection
 & \rr PyWake decides every acceptance
 & \rr PyWake calls \tabularnewline
\bottomrule
\end{tabular}
\end{table}

\paragraph{Site selection on an interference matrix.}
Binary site-selection models were initiated by Donovan~\cite{donovan2005}
and Fagerfj\"all~\cite{fagerfjall2010}, and systematized by Turner et
al.~\cite{turner2014}. They price the wakes through pairwise deficits
computed before the optimization (or between solves, as in Kuo et
al.~\cite{kuo2016}), and so does the neighborhood-search MILP of
P\'erez-R\'ua et al.~\cite{perezrua2023}. Fischetti and
co-authors~\cite{fischetti2016proximity,fischetti2021nextgen,%
fischetti2023integrated} carried the line into industrial
practice~\cite{fischetti2020vattenfall}, the field is surveyed
in~\cite{fischetti2019overview}, and EndoWake is introduced
in~\cite{companion}. Cazzaro and Pisinger~\cite{cazzaro2022vns} place
hundreds of turbines among thousands of sites by variable neighborhood
search, and Cazzaro et al.~\cite{cazzaro2022multiscale} also choose the
shape of the farm; Thomas et al.~\cite{thomas2023} compare discrete and
continuous methods on one problem, and Baker et al.~\cite{baker2019}
benchmark continuous ones. The family expresses counting and logical
conditions linear in the binaries, while the wake is frozen as data.

\paragraph{Continuous coordinates.} TOPFARM~\cite{topfarm}, built on
PyWake, and FLORIS~\cite{floris} move the turbine coordinates along
gradients of annual energy, in the line of the exact wake gradients of
Guirguis et al.~\cite{guirguis2016} and of the boundary-grid
parameterization of Stanley and Ning~\cite{stanley2019}. Faster variants
parallelize on large farms~\cite{valottarodrigues2024speeding}, sample the
wind conditions at random~\cite{quick2023sgd}, or integrate the rose by
quadrature~\cite{padron2019pce}. FLOWERS~\cite{locascio2022,locascio2024}
expands the rose in a Fourier series, so that annual energy becomes a
closed-form, differentiable function of the positions~\cite{locascio2025}:
it removes the rose discretization and keeps linear superposition and
free-stream thrust, whereas EndoWake keeps a discretized rose, propagates
from the effective speed at every machine and implies, rather than
prescribes, a linear rule only where two wakes merge sideways.

\paragraph{CFD in the loop.} King et al.~\cite{king2017} and Antonini et
al.~\cite{antonini2020adjoint} optimize against a flow field carried as
the state of an adjoint computation, the latter with a RANS on a farm in
complex terrain. Gnegel et al.~\cite{gnegel2021} find that writing a
discretized linear partial differential equation into a mixed-integer
program does not scale. EndoWake writes a discretized propagation into the
program too, but a local, one-sided one that a forward sweep solves, and
calls the heavy model once, offline, to calibrate it
(Section~\ref{sec:distill}).

\paragraph{Site constraints.} The award of an area commonly comes with
conditions of four kinds~\cite{fischetti2019overview,fischetti2020vattenfall}:
bounds on the number of turbines, fixed costs per site charged only if
the site is used, non-convex exclusions and sheltering structures, and a
cardinality left to the optimization. Costs and exclusions are written
exactly on continuous coordinates~\cite{criadorisco2024zones}; the counts
ask for a discrete treatment, which in TOPFARM comes through a custom
component or a penalty~\cite{topfarm}. On binary
site variables all four are linear constraints, as on an interference
matrix; what EndoWake adds is the field in the same variables, so that
counts, site costs, obstacles and wakes sit in one formulation
(Section~\ref{sec:mip}, tested in Section~\ref{sec:exp-sidecon}).

\paragraph{Cheaper models of a heavier one.}
Fuga~\cite{ott2011fuga,ott2014fuga} linearizes the steady RANS equations
about an undisturbed boundary-layer profile and reconstructs a farm from
single-turbine look-up tables~\cite{pena2018anholt}.
TurbOPark~\cite{nygaard2020turbopark,pedersen2022turbopark} lets the wake
expansion follow the local turbulence intensity, a mechanism that a fixed
kernel calibrated on one isolated wake cannot carry
(Section~\ref{sec:mermaid}). Neural, graph and Gaussian-process surrogates
learn a wake, or the power of a farm, from
RANS~\cite{ti2020mlwake,li2023graphwake,li2022sfnet}, large-eddy
simulation~\cite{bempedelis2024bo} or an engineering
model~\cite{park2019pgnn,anagnostopoulos2023wakenet}, a
physics-informed network of the wake drives a differential-evolution
layout search in Li et al.~\cite{li2025hybrid}, a Kriging surrogate of
the CFD-computed energy of sampled layouts drives a genetic layout search
in Wang et al.~\cite{wang2024kriging}, and Shao et al.~\cite{shao2025pce}
cut the cost of the annual energy by polynomial chaos and search the
layout by Kriging with expected improvement. All of them deliver
an evaluator: a local search on it carries the constraints above as
feasibility checks, but none supplies the bounds on the objective that a
MILP formulation does (Section~\ref{sec:mip}). The
diffusion-based wake model of Ali et al.~\cite{ali2024diffusion} reads a
single wake as the diffusion of a passive scalar from a disk, with a
streamwise scaling that conserves momentum; the kernel of
Section~\ref{sec:notation} is a discrete diffusion too, but fitted, and
conserving nothing.

Multifidelity optimization is the closest precedent of
Section~\ref{sec:proxy}: TOPFARM itself began as a multi-fidelity
framework, in which cheaper models of wakes and costs take the first
iterations and costlier ones refine the layout~\cite{rethore2014topfarm}.
On a layout case, Jasa et
al.~\cite{jasa2022multifidelity} let a Jensen model, corrected by a
surrogate of its discrepancy, steer a gradient-based search within a trust
region, and call a Gauss--curl hybrid model to update the correction; they
find an optimum close to the high-fidelity one, not guaranteed to coincide
with it. Zhang et al.~\cite{zhang2026} report the other side of the
trade-off: at a fixed time budget, a costlier wake model of similar
accuracy allows fewer evaluations and yields worse layouts. Guiding
a search with a cheap model toward the optimum of a finer one is older:
Kaja Kamaludeen et al.~\cite{kajakamaludeen2016} map a coarse wake model
onto a fine one by manifold mapping and reach the optimum of the fine
model at a fraction of its evaluations, and Song et
al.~\cite{song2023surrogate} pre-screen with a learned surrogate inside a
population search, for the control of a floating farm rather than its
layout. Closest to the directional error of Section~\ref{sec:sirens},
Kim~\cite{kim2026gnn} learns, with a graph neural network, the
correction from a coarse evaluation on sixteen directions to the fine one
on $360$, and applies it inside the layout optimization: the approximate
objective is corrected, whereas here every accepted move is confirmed on
the full one. What Section~\ref{sec:proxy} adds is a specific combination: a linear
propagator distilled from the single wake of PyWake ranks whole
neighborhoods in batch on a GPU, the shortlist is sized from measured
evaluation rates, and PyWake decides every acceptance. Apart from
PyWake and TOPFARM, the tools cited in this section are positioned here
and not run, and we claim nothing about their running time. 
\section{EndoWake, and the setup of the experiments}
\label{sec:surrogate}

In this section we recall EndoWake (Sections~3 and~4 of~\cite{companion})
at the depth needed to follow the experiments, and we fix the setup shared
by all later sections. The wind speed of every cell is a decision variable,
carried downstream by a linear stencil with three parameters, so a
planner's constraints enter the model without touching the physics. One
forward sweep per wind sector computes the fastest field the model admits,
and this sweep is the evaluator behind every EndoWake search of the paper.

\subsection{Wind as a variable, and the stencil that moves it}
\label{sec:notation}

Consider a regular $N \times M$ grid with spacing $\Delta$, expressed in
rotor diameters $D$ (in our experiments $\Delta = D$, $D/2$ or $D/3$). To
ease exposition, we first let the wind blow from the north with undisturbed
speed $\wmax$, so that the cell immediately downstream of $(i,j)$ is
$(i+1,j)$. A turbine removes a fraction of the speed that reaches it, and
the resulting \emph{deficit} $d = \wmax - w$ travels downstream one plane of
cells at a time, through a fixed \emph{stencil}. The rule mixes the speed of
the cell upwind with those of the adjacent and the second columns, and adds
a recovery term that pulls the field back toward the free stream:

\begin{equation}
w_{ij} \;=\; \gamma \wmax + (1-\gamma)
\sum_{m=-2}^{2} K_m \, w_{i-1,j+m},
\qquad
K = \Bigl(\tfrac{\beta}{2},\tfrac{\alpha}{2},1-\alpha-\beta,
          \tfrac{\alpha}{2},\tfrac{\beta}{2}\Bigr).
\label{eq:prop}
\end{equation}

Here $\alpha,\beta \ge 0$, with $\alpha+\beta \le 1$, control lateral
diffusion to the adjacent and to the second columns, and $\gamma \in [0,1]$
controls recovery. The coefficient vector $K$, the \emph{kernel}, is the
same on every cell and sums to one, so the right-hand side
of~\eqref{eq:prop} is a convex combination of the free-stream speed and of
five speeds upwind: the rule can neither create wind nor drive it
negative. An index $w_{i-1,j+m}$ that leaves the grid stands for $\wmax$.
Evaluating a layout is one sweep over the grid, five multiply--add
operations per cell, which takes milliseconds on a laptop at industrial
size (Section~\ref{sec:exp-big}).

\subsection{The model}
\label{sec:mip}

We next embed the propagation rule into the MILP model of Section~4
of~\cite{companion}. Let $\mathcal{C}$ be the set of grid cells and
$\mathcal{A} \subseteq \mathcal{C}$ the allocatable cells, a proper subset
of $\mathcal{C}$ under the non-convex geometries of
Section~\ref{sec:exp-sidecon} and wherever obstacles are present. Let
$\mathcal{S}$ be a set of wind scenarios, scenario $s$ having direction
$\theta_s$, undisturbed speed $\wmax_s$ and probability $p_s$; let $C_c$ be
the cost of a turbine at $c$, in units of expected power, and $P_{\max}$
the rated power. Finally, let $\rho \in [0,1)$ be the fraction of the speed
that survives on the footprint cells of a rotor: rather than setting it to
its momentum-theory value, we fit the source amplitude $1-\rho$ together
with the stencil, at every spacing (Section~\ref{sec:calibration}). The
variables are

\[
\begin{array}{ll}
x_c \in \{0,1\}, & c \in \mathcal{A}
   \quad\text{(a turbine is installed at $c$)},\\[2pt]
w^s_c \in [0,\wmax_s], & c \in \mathcal{C},\ s \in \mathcal{S}
   \quad\text{(wind speed)},\\[2pt]
z^s_c \in [0,P_{\max}], & c \in \mathcal{A},\ s \in \mathcal{S}
   \quad\text{(power output)}.
\end{array}
\]

For every scenario $s$ and cell $c=(i,j)$, let $\mathrm{up}_s(c,m)$,
$m = -2,\dots,2$, be the five upwind neighbors of $c$ in the stencil of
$s$. The \emph{footprint} $\mathcal{F}_s(c)$ of a turbine at $c$ is the set
of $2h+1$ cells of the next plane centered downstream of $c$, i.e.,
$\{(i+1,\,j+\ell) : |\ell| \le h\}$ for a wind from the north, where the
half-width $h$ is $R/\Delta - \tfrac12$ rounded to the nearest integer,
ties to even, and $R = D/2$ is the rotor radius. The footprint spans the
rotor exactly when $D/\Delta$ is odd; at $\Delta = D/2$ it is the single
cell directly downstream ($h = 0$). Let $(v_b, P_b)$, $b = 1, \dots, B$, be
the breakpoints of the power curve, $\widehat{P}$ the piecewise-linear
interpolant they induce and $\lambda^s_{cb} \ge 0$ the interpolation
weights. Finally, for a minimum separation $\delta$ (say), let $N(c)$
collect the allocatable sites closer than $\delta$ to $c$. Our model reads:

\begin{equation}
\max \ \sum_{s \in \mathcal{S}} p_s \sum_{c \in \mathcal{A}} z^s_c
\;-\; \sum_{c \in \mathcal{A}} C_c x_c .
\label{eq:mip-obj}
\end{equation}

\begin{equation}
w^s_{c} \;\le\; \gamma_s \wmax_s + (1-\gamma_s)
\sum_{m} K^{s}_m\, w^s_{\mathrm{up}_s(c,m)} .
\label{eq:mip-flow}
\end{equation}

\begin{equation}
w^s_{c'} \;\le\; \rho\, w^s_{c} + \wmax_s\,(1-x_c),
\qquad c \in \mathcal{A},\ c' \in \mathcal{F}_s(c).
\label{eq:mip-wake}
\end{equation}

\begin{equation}
\begin{aligned}
&\sum_{b=1}^{B} \lambda^s_{cb} \;=\; 1, \qquad
\sum_{b=1}^{B} v_b\, \lambda^s_{cb} \;=\; w^s_c, \qquad
(\lambda^s_{c1}, \dots, \lambda^s_{cB}) \ \text{SOS2}, \\[2pt]
&z^s_c \;\le\; \sum_{b=1}^{B} P_b\, \lambda^s_{cb}, \qquad
z^s_c \;\le\; P_{\max}\, x_c .
\end{aligned}
\label{eq:mip-power}
\end{equation}

\begin{equation}
\sum_{v \in N(c)} x_v \;\le\; |N(c)|\,(1-x_c) \quad \forall c \in \mathcal{A},
\qquad
N_{\min} \le \sum_{c \in \mathcal{A}} x_c \le N_{\max}.
\label{eq:mip-dist}
\end{equation}

The objective~\eqref{eq:mip-obj} maximizes the expected power over the
scenarios, net of the installation costs. Constraints~\eqref{eq:mip-flow}
propagate the field plane by plane in every scenario, an index off the grid
reading $\wmax_s$; $K^{s}$ and $\gamma_s$ are the kernel and the recovery
rescaled to the direction $\theta_s$ (Section~\ref{sec:rose}).
Constraints~\eqref{eq:mip-wake} are the source of the wake: a turbine at
$c$ caps the speed on each footprint cell at $\rho$ times the speed that
reaches $c$, while the big-$M$ term $\wmax_s(1-x_c)$ deactivates the
constraint when $x_c = 0$. Note that a turbine scales down whatever speed
reaches it, instead of resetting it to a fixed value: a machine deep in a
wake casts a proportionally weaker wake, and a line of machines needs no
superposition rule. A rule is implied, not chosen, where two wakes merge
sideways: the propagation is linear, hence the deficits of two neighbors
add on the cells they share, and~\eqref{eq:mip-wake} caps the sum at the
next rotor (Section~\ref{sec:exp-fidelity} measures what this gives on two
merging wakes). Block~\eqref{eq:mip-power} converts speed into power: under
the \emph{SOS2} condition~\cite{beale1970} (at most two positive weights,
and consecutive if two), it states $z^s_c \le \widehat{P}(w^s_c)$, with no
binary variables added. The interpolant has $B = 8$ breakpoints in all
experiments, with the jump at the cut-in speed of $3$\,m/s reproduced
through two breakpoints $10^{-3}$\,m/s apart. It enters the MILP
only: the local searches, Hexaly and the batch evaluator of
Section~\ref{sec:exp-batch} score with the exact curve, and the layout
returned by the SOS2 model, re-evaluated with the exact curve, is within
1.13\%\ of the best layout known for the same
instance (\TN). Finally,
constraints~\eqref{eq:mip-dist} impose the minimum separation and bound
the number of machines; as in~\cite{companion}, our implementation
replaces the aggregated form by the equivalent clique inequalities on the
conflict graph.

\label{sec:fixpoint}%
Note that~\eqref{eq:mip-flow}, \eqref{eq:mip-wake}
and~\eqref{eq:mip-power} are inequalities, so every block is a relaxation
in the direction the objective pushes. Proposition~1 of~\cite{companion}
proves that this costs nothing: the fields
satisfying~\eqref{eq:mip-flow} and~\eqref{eq:mip-wake} have a greatest
element $w^{*}(x)$, computed in $O(|\mathcal{C}|)$ per scenario by one
forward sweep in topological order, and every objective non-decreasing in
$w$ attains its maximum there. The sweep of Section~\ref{sec:notation} is
therefore the exact evaluator of the model, for the local search and the
Hexaly model~\cite{hexaly} of Section~\ref{sec:exp-heur} as for the batch
GPU evaluator of Section~\ref{sec:exp-batch}. The exception is a variant of
the field, the \emph{one-cell correction}: where the footprint is a single
cell ($\Delta = D$ and $D/2$), arm~B of Section~\ref{sec:rans-test} and
the runs of Section~\ref{sec:cones} read the speed of a turbine at its true
position on an oblique flow line and split its cap between the two cells
that bracket the point where that line crosses the next plane (\TN\
specifies it). Every other search, arm~B$'$ at $\Delta = D/3$ included,
uses the footprint of~\eqref{eq:mip-wake}.

Overall, the model has $|\mathcal{A}|$ binary variables and
$O(|\mathcal{S}|\,|\mathcal{C}|)$ continuous variables and constraints
besides the separation family, against the $\Theta(|\mathcal{A}|^2)$
interaction terms of interference-matrix formulations. The physics lives
entirely in~\eqref{eq:mip-flow}--\eqref{eq:mip-wake}: any condition a
planner can write on $x$, such as a cardinality per zone, a non-convex
exclusion, a per-site cost or an integer count of machines of a given
class, joins~\eqref{eq:mip-dist} without touching the field
(Section~\ref{sec:opt-tool}), and a solid obstacle is a cell with $w^s_c =
0$ as an upper bound, whose shadow the propagation produces.

\subsection{Wind roses}\label{sec:rose}
Each sector of a measured rose is a scenario with its own direction, mean
speed and frequency. For an oblique direction, propagation advances plane
by plane along the \emph{dominant axis}, the coordinate with the larger
component of the wind vector, so that every cell depends only on the
previous plane and a scenario is evaluated by one sweep, as an aligned
wind is. An integer staircase follows the lateral drift, and the kernel is
rescaled so that diffusion and recovery per unit length match the aligned
case: with $\mu = 1/|u_\parallel| \in [1,\sqrt2]$ the length of a plane
step in cells, $u_\parallel$ being the component of the unit flow vector
along the dominant axis, $\alpha$ and $\beta$ become $\mu^3\alpha$ and
$\mu^3\beta$, and $1-\gamma$ becomes $(1-\gamma)^{\mu}$ (\TN\ gives the
algebra).

\subsection{Experimental setup}
\label{sec:exp}

Four words recur in what follows. The \emph{evaluator} is the code that
scores a layout under a given model. A \emph{guide} is the model a search
calls while it moves machines. A \emph{judge} scores each layout once,
after the search that produced it has ended. An \emph{arm} is one search
pipeline in a comparison: the arms start from the same layouts and share
the power curve, the calibration data and the judge, and differ in the
guide (the wake model and the way it integrates the wind climate) and,
where a section says so, in the scan of the neighborhood and in the
hardware; a comparison of arms measures complete pipelines, not the wake
model in isolation.

The turbine is that of~\cite{companion}: a machine of the 3\,MW class with
rotor diameter $D=112$\,m, cut-in $3$\,m/s and cut-out $25$\,m/s. To avoid
proprietary power-curve data, we give it the idealized curve $P(v) =
\min\{\tfrac12 \rho_{air} A C_p v^3, P_{\max}\}$ with $\rho_{air} =
1.225$\,kg/m$^3$ and $C_p = 0.45$, which rates at
$10.34$\,m/s. Four exceptions run instead the
SWT-2.3-93 of Lillgrund ($D = 93$\,m, hub
height 65\,m): (i)~the RANS reference of
Section~\ref{sec:rans}, at its own inflow speed; (ii)~the arms of
Section~\ref{sec:rans-test}; (iii)~PyWake where it judges Lillgrund
(Section~\ref{sec:sirens-judge}); (iv)~both models in the Lillgrund experiments
of Section~\ref{sec:proxy}. The Horns Rev~1 experiment of the same section
uses the V80 curves supplied with PyWake. Every other search on Lillgrund
keeps the idealized curve, with its geometry measured in the diameter of
that site.

The grid spacing is $\Delta = D$ wherever the RANS is the judge
(Section~\ref{sec:sirens}, which also runs $D/3$ and says where) and
$\Delta = D/2$ in the optimization experiments of
Section~\ref{sec:opt-tool}; every table states which. The minimum
separation is $3\,D$, and the undisturbed speed is $10$\,m/s in every
sector of a synthetic rose, whose sectors differ only in direction and
weight. Installation costs are zero ($C_c = 0$). We depart from these
choices only where the departure is itself the experiment: from the zero
costs in the free-cardinality study of Section~\ref{sec:exp-sidecon} and in
the lower panel of Table~\ref{tab:biggrids}, and from the $3\,D$ separation,
which becomes the farm-like $5\,D$\ in
Section~\ref{sec:exp-big}. The measured climates of Horns~Rev~1 (\TN) and
Lillgrund (Section~\ref{sec:exp-lillgrund}) replace the synthetic wind,
each sector at its own Weibull mean speed and frequency. Where the
reference is analytical, the stencil parameters come from a calibration
against the Gaussian model of Bastankhah and
Port\'e-Agel~\cite{bastankhah2014} as implemented in
PyWake~2.6.20~\cite{pywake}, tabulated in \TN; the RANS runs and their
distillate are the subject of Section~\ref{sec:distill}.
Appendix~\ref{app:exp-conventions} says which machine carried what, with
the software versions and the numerical conventions; time limits, thread
counts and seeds are stated in the caption of each table. Except for arm~A
of Section~\ref{sec:sirens} and the annual energy of
Section~\ref{sec:proxy}, which integrate the power over the Weibull
distribution of each sector, all energies are expected power in kW
averaged over the scenarios of the rose, each at the speed it carries. An
annual energy is $8\,760$ hours times the corresponding expected power, so
relative gains are the same in both scales. 
\section{How a heavy CFD model is distilled into EndoWake}
\label{sec:distill}

In this section we distill a steady RANS computation into the four numbers
of EndoWake, fitted on the wake of one isolated turbine, and we measure the
distillate in fit, in speed and in fidelity on interacting machines. The
RANS reference built here is also the judge of Section~\ref{sec:sirens}.

\subsection{The recipe: what a calibration consumes}
\label{sec:calibration}

A steady RANS computation of the whole farm is the most faithful model we
can afford for scoring final layouts but, at minutes per wind direction on
a compute node (Table~\ref{tab:rans_setup}), not for the millions of a
search. By \emph{distilling} it we mean something narrow: we calibrate the
four numbers of EndoWake, the stencil weights $(\alpha,\beta,\gamma)$ and
the source amplitude $1-\rho$ of Section~\ref{sec:mip}, once, and the
propagation then generates the interactions of a whole farm at the cost of
a linear sweep. The calibration consumes one object, the deficit field
$d^{\mathrm{ref}}$ of a single isolated turbine, sampled on a window from
$3\,D$\ to $12\,D$\
behind the rotor, which covers the spacings an optimizer produces and
leaves out the near wake, and minimizes a weighted least-squares
discrepancy between the propagated single-turbine field and
$d^{\mathrm{ref}}$ over the window; the residual is the root-mean-square
error relative to the reference deficit, with the protocol
of~\cite{companion} kept unchanged at every cell size.

The fit is cheap because the field is linear in the source amplitude, and
the deficit $d = \wmax - w$ injected as a transverse profile $S$ at the
plane of the footprint obeys, $\ell$ planes downstream,
\begin{equation}
d_{\,i_0+\ell,\;\cdot} \;=\; (1-\gamma)^{\ell} \; \bigl(S * K^{*\ell}\bigr),
\label{eq:closedform}
\end{equation}
where $*$ is the discrete convolution along the plane and $K^{*\ell}
= K * K * \cdots * K$ ($\ell$ factors) the $\ell$-fold convolution power of
the stencil, as long as the wake has not reached the lateral edges of the
grid (\TN\ gives the proof). At every trial $(\alpha,\beta,\gamma)$ the
best amplitude is thus available in closed form, and the calibration is a
three-parameter search that a derivative-free method settles in seconds.
The closed form serves the fit only: in a farm the cap
of~\eqref{eq:mip-wake} is applied after every plane, a nonlinear
operation, so no single convolution carries the whole farm; the sweep is
a direct five-tap convolution per plane, the cheapest for this stencil.

The amplitude $1-\rho$ is not the axial induction factor $a$ of
momentum theory: it is a lumped strength for the source the model
injects, a row of $2h+1$ cells on the plane behind the rotor, from which
the recovery already acts, fitted on the window from three diameters on.
Its value depends on the cell, since the width of the source does
(Section~\ref{sec:rans-calib}), and amplitudes at different cells are
comparable neither with one another nor with $2a$. Holding the amplitude
at the momentum value costs, at $\Delta = D$,
8.0\%\ of field error against
7.7\%\ on the analytical target of
Section~\ref{sec:exp}, the Gaussian wake of~\cite{bastankhah2014} at $k =
0.04$, on which the recipe leaves 7.7\%\ at
$\Delta = D$ and 45\%\ at the $\Delta = D/2$
of the optimization experiments, where the source is a single cell half a
rotor wide (\TN): the comparisons that use this kernel remain legitimate,
since every solver evaluates the same implemented field, but they do not
inherit the accuracy of the fit at $\Delta = D$. The calibrated parameters
depend on the cell size, so a layout obtained at one spacing must be
re-evaluated, never re-interpreted, at another. Two discrepancies are
structural whatever the parameters: a kernel of finite support widens as
the square root of the downstream distance, where wakes widen linearly,
and the recovery $\gamma$ decays the deficit exponentially, where wakes
decay polynomially. And a kernel fitted on one isolated wake carries
none of the physics at the scale of a farm: blockage, shear and veer,
meandering and added turbulence are absent, and what the distillate gets
right on interacting machines must be measured
(Section~\ref{sec:exp-fidelity}).

\subsection{The reference, and the fit on its single wake}
\label{sec:rans}
\label{sec:rans-oracle}

Table~\ref{tab:rans_setup} states the setup of the steady RANS of the
Lillgrund farm that we call \emph{the RANS reference}, and never ``the
CFD'': it is one closure on one mesh, not the discipline.
Appendix~\ref{app:rans-details} gives the verification of the closure, the
average over the wind direction, the monitor, the grid convergence and the
run-to-run noise, and \TN\ the record of the runs.

\begin{table}[htbp]
\centering\small
\caption{The RANS reference in eight rows. It is the judge of
Section~\ref{sec:sirens}.}
\label{tab:rans_setup}
\begin{tabular}{l p{0.78\textwidth}}
\toprule
solver & OpenFOAM v2412~\cite{openfoam}, steady incompressible RANS
(SIMPLE), up to $1\,500$ iterations per direction. The single-wake and
calibration runs of Section~\ref{sec:rans} stop at the residual control
$p < 10^{-4}$, $U$, $k$, $\varepsilon < 10^{-5}$, reached after
505\ to
604 iterations; every judged layout
runs at a tenfold tighter control
($p < 10^{-5}$, $U$, $k$, $\varepsilon < 10^{-6}$), which
24\%\ of the
552\ judged runs of the judgment
campaign of this paper (all its arms, seeds and directions; the
240%
\ runs of the controlled comparison
of Section~\ref{sec:sirens-judge} are audited separately in \TN)
reach within the $1\,500$ iterations
(Appendix~\ref{app:rans-details}) \\
domain & $8\,D$ upstream, $15\,D$ downstream, $6\,D$ to either side and
$800$~m in height, with the wind along the positive $x$ axis \\
boundaries & inlet: the atmospheric-boundary-layer profiles of the solver
library for $U$, $k$ and $\varepsilon$, with the constants of the
closure; outlet: fixed pressure, \texttt{inletOutlet} on the
transported fields; sides and top: slip; ground: no slip with a rough-wall
function at $z_0$, whose $\varepsilon$ condition keeps the default
constants of the library ($C_\mu = 0.09$, $\kappa = 0.41$); along the
empty domain the hub-height speed over the farm drifts by at most
1.2\%\ (Appendix~\ref{app:rans-details}) \\
mesh & $25$~m hexahedral cells, refined to $12.5$~m around the farm
(about seven cells per rotor diameter) \\
inflow & neutral logarithmic profile over sea, $z_0 = 2\cdot 10^{-4}$~m,
$8$\,m/s at hub height, one run per direction; the inlet $k$ and
$\varepsilon$ follow the log law with $k = u_*^2/\sqrt{C_\mu}$, which puts
the turbulence intensity at hub height at
6.2\% \\
turbines & actuator disks with the SWT-2.3-93 thrust and power curves,
read against the mean velocity of nine points over the rotor one
diameter upstream \\
closure & $k$-$\varepsilon$-$f_P$ of van der Laan et
al.~\cite{vanderlaan2015kefp,vanderlaan2015double,vanderlaan2015farms}, implemented by us
from the published formulation and verified on the single-turbine
case of~\cite{baungaard2022}, distributed with~\cite{actuatordiskfoam};
$C_\mu = 0.03$, $C_1 = 1.21$, $\kappa = 0.40$, $\sigma_\varepsilon =
1.30$ and $C_R = 4.5$, the constants of the authors: none is tuned \\
cost & 405~s per direction on 55 cores
of one node, mesh included \\
\bottomrule
\end{tabular}
\end{table}

The turbulence closure decides the numbers. On the aligned row at
$222^\circ$ (spacing $4.3\,D$), standard $k$-$\varepsilon$ with
atmospheric constants recovers the wake too fast: the second turbine
produces 0.44\ of the front-row power on
our grid at that single direction, a value that the averaging over
directions discussed below can only raise. The $k$-$\varepsilon$-$f_P$ closure delays the recovery by
scaling the eddy viscosity with the local shear; its constants are those
of its authors (Table~\ref{tab:rans_setup}). On the single-turbine case
of~\cite{baungaard2022}, distributed with~\cite{actuatordiskfoam},
our implementation follows the field of the authors' solver to
0.033\ of the free-stream speed along the
wake axis, in root mean square between $3$ and $10\,D$
(Appendix~\ref{app:rans-details}). The previous submission of this
work used $C_R = 60$ to compensate for an $f_P$ factor left out of the
production term of the $\varepsilon$ equation; we corrected the
implementation and recomputed the calibration, the judgments and the
controls with $C_R = 4.5$ (Appendix~\ref{app:rans-details}).

Against the SCADA of the site the comparison needs care, because a
measured value is an average over a bin of wind directions, and the wind
direction is itself uncertain~\cite{gaumond2014}. At $222^\circ$ the
SCADA, binned over $\pm 2.5^\circ$ at $8 \pm 0.5$\,m/s, put the second
turbine at about 0.43\ of the front-row
power and the machines from the fifth on at
0.30--0.35~\cite{sebastiani2021}. At the
single direction, the bottom of a narrow dip, the reference gives
0.23\ at the second machine and settles
down the row at 0.33\unskip; averaged
over the direction with a Gaussian weight of standard deviation
$3.3^\circ$, the value estimated for this farm
in~\cite{vanderlaan2015farms}, it gives 0.31\
at the second machine and 0.38\ from the
fifth on. At $120^\circ$, the other aligned direction and a tighter one
($3.3\,D$), it gives 0.17\ at the single
direction and 0.21\ averaged, against a
measured value of about 0.2\unskip. What
remains is largest at the second machine of the $4.3\,D$ row, where the
reference gives less power than the site records: as a judge it
penalizes, if anything, a closely aligned pair. At more than an hour of a
compute node per layout at twelve directions, it can judge a search but
never sit inside one.

\paragraph{One field, two models.}
\label{sec:rans-calib}
We now run the recipe on the RANS field of an isolated turbine, with the
footprint the model injects, and for an apples-to-apples comparison we fit
the Gaussian model on the same data (Table~\ref{tab:rans_calib}), with
its one free parameter, the wake-expansion coefficient, fitted and at the
$k = 0.04$ of the analytical target.

\begin{table}[t]
\centering
\caption{Calibration on the RANS deficit field of an isolated turbine ($k$-$\varepsilon$-$f_P$, $C_R = 4.5$, 8\,m/s, window 3--12\,$D$) at four cell sizes $\Delta$: $\alpha$ and $\beta$ are the weights of the adjacent and second columns ($\beta = 0$ at every optimum), $\gamma$ the decay per plane, amplitude the deficit injected at the source; the Gaussian wake of Bastankhah and Port\'e-Agel is fitted on the same field through its growth rate $k$ (column $\gamma$), with the fitted $k^*$ and with the literature value. Residual: root-mean-square error of the deficit over the window, relative to the RANS deficit.}
\label{tab:rans_calib}
\begin{tabular}{llrrrrr}
\toprule
model & $\Delta$ & $\alpha$ & $\beta$ & $\gamma$ / $k$ & amplitude & residual \\
\midrule
linear kernel & $D$ & 0.0328 & 0.00 & 0.1035 & 0.644 & 9.4\% \\
linear kernel & $D/2$ & 0.1094 & 0.00 & 0.0291 & 1.000 & 15.6\% \\
linear kernel & $D/3$ & 0.0783 & 0.00 & 0.0292 & 0.579 & 6.7\% \\
linear kernel & $D/4$ & 0.0629 & 0.00 & 0.0242 & 0.491 & 11.3\% \\
\midrule
Gaussian, fitted $k^*$ & --- & --- & --- & 0.0337 & --- & 7.3\% \\
Gaussian, $k = 0.04$ & --- & --- & --- & 0.04 & --- & 15.5\% \\
\bottomrule
\end{tabular}
\end{table} 

Our field is the closer of the two at $\Delta = D/3$ only, with the
obvious qualification that our fit has four free numbers and theirs one;
at $D$, $D/4$ and $D/2$ the fitted Gaussian model is the more faithful.
On a single wake the
distillate is thus about as faithful as a fitted Gaussian model, not more;
what sets it apart is the linear form of Section~\ref{sec:mip} and the
throughput of Section~\ref{sec:exp-batch}. The poor fit at $\Delta = D/2$
is geometric: the source is a row of cells centered on the machine, whose
width matches the rotor only when an odd number of cells covers one
diameter ($\Delta = D, D/3, D/5, \dots$); at $D/2$ it is a single cell
half a rotor wide, at $D/4$ five cells a quarter too wide, and no choice
of parameters corrects either.

The weight $\beta$ vanishes in every fit of Table~\ref{tab:rans_calib}:
over this window the second neighbor does not help, and forcing
$\beta = 0.01$ or $0.05$ at $\Delta = D$ raises the residual from
9.4\%\ to 10.5\%\
and 23.4\%\unskip. We keep $\beta$ in the
model because the kernel is the general five-point stencil
of~\cite{companion}, on which a reference with a wider or flatter wake
would load it.

\paragraph{Sensitivity to the inflow.} Refitted on single-wake fields of
the same reference over two rougher surfaces ($z_0 = 0.005$ and $0.05$~m,
hub-height turbulence intensities 8.3\%\ and
10.9\%\ against
6.2\%\unskip) and at $6$ and $10$\,m/s (\TN),
the kernel is insensitive to the speed and sensitive to the turbulence: at
the roughest surface $\alpha$ grows by a factor
1.7\ and the amplitude falls from
0.644\ to 0.483\unskip,
with a residual at $\Delta = D$ of at most 9.4\%\
in every condition. The kernel must thus be fitted at the turbulence of
the site. In decisions this costs less than in fit: under the
six-direction judge of Section~\ref{sec:cones}, arm~B of
Section~\ref{sec:rans-test} gains +3.70\%\ with the
kernel of the roughest surface and +4.01\%\ with the
kernel of the site, on the same 2\ seeds.
The farm reference and every judgment ran at one stability and one
speed.

\subsection{The distillate on a GPU}
\label{sec:exp-batch}

The maximal-field property of Section~\ref{sec:mip} makes each evaluation
a single sweep, computed plane by plane in linear time: the same
five-point convolution on every plane, followed by a cap on the cells
downwind of a turbine, with no data-dependent branching, so $n_b$
candidate layouts are scored by one sweep on arrays of shape $n_b\times M$
in place of $M$. On one node of the university cluster, on a
100$\times$100\
grid carrying 30\ turbines, batching alone is
worth a factor 26\ on one CPU core in
double precision, and the GPU of the node, an
NVIDIA RTX PRO 6000 Blackwell Server Edition\unskip, adds in single precision a further
47.9\unskip$\times$, for
1\,093\,103\ layouts per second on one direction (\TN).
A data-center card is not a requirement: under the twelve sectors of
Lillgrund, the entry-level mobile GPU of the laptop of
Appendix~\ref{app:exp-conventions} sustains 60\,130\
layouts per second with the sweep of a plane fused into a single kernel.
These figures are batch
throughputs that differ in rose, card and grid, and none is the speedup of
a search, which Section~\ref{sec:proxy} measures.
At equal hardware the distillate is not faster than PyWake. On one
node of the cluster of Section~\ref{sec:proxy} and on the instance of
that section (4\,728\ sites, $360$ directions,
one speed), the vectorized CPU form of the evaluator scores, in double
precision, 27%
\ layouts per second
on one core against 58%
\ for PyWake, and
1152%
\ against
1317%
\ on thirty processes. The advantage
is the GPU of the same node, which scores
31196%
\ layouts per second in single
precision and 17130%
\ in double,
24%
\ times PyWake on thirty processes:
the speedups of Section~\ref{sec:proxy} are those of a GPU added to the
processes of PyWake, not of one evaluator against the other at equal
resources.

Single precision is the natural arithmetic for EndoWake: its rounding,
about $10^{-7}$, lies far below the error of the model itself
(9.4\% against the RANS at $\Delta = D$);
the energies move by $2\times10^{-6}$\ in relative terms,
and 0\ of the
496\ pairwise comparisons we checked change
sign, while in double precision the batched sweep reproduces the
reference evaluator to $2\times10^{-16}$\unskip. Note
that the distillate is an accelerated RANS only in the sense that it was
distilled from one RANS wake: it returns the hub-height field on one
plane, for a discrete rose of directions at one speed, with no vertical
structure, no turbulence quantities and no transient.

\subsection{What the distillate gets right on interacting machines}
\label{sec:exp-fidelity}
\label{sec:exp-pywake}

Fitted on one isolated wake, the distillate must still be right on
interacting machines. Three tests on configurations it was not fitted
to, with tables, protocols, seeds and sweeps in \TN, bound what it gets
right. (i)~\emph{Two wakes that merge.} On two turbines abreast at
lateral offsets of $1$, $1.5$ and $2\,D$, our kernel overstates the
transverse integral of the merged deficit at every station behind them
and at both cell sizes, by 1.04\ to
1.21\ times the RANS integral, and its
largest profile error (5.6\%\ of the
free stream at $\Delta = D/3$, 9.3\%\ at
$\Delta = D$) is above those of the Gaussian model summed linearly or in
root-sum-square (4.2\%\ and
3.5\%\unskip).
(ii)~\emph{Higher-order interaction along a row}, with PyWake held fixed
as the reference: a row of aligned turbines rebuilt from pairwise
coefficients measured on PyWake with two turbines in free stream is in
error by 95.9\%\ at
3\,D\unskip, where EndoWake is in error
by 5.8\%\unskip; the advantage is not
monotone in spacing. (iii)~\emph{The ordering of local search moves.} On
pairs of layouts that differ by one relocated turbine, judged by PyWake
under the twelve-sector rose of Horns Rev, the endogenous field gets the
sign of a move right more often than the Jensen pairwise models and less
often than the pairwise models whose coefficients are measured on the
judge, and on Kendall's $\tau_b$ over all the pairs it is the lowest of
the five evaluators; the pairs come from a few searches per condition and
are not independent replicates. 
\section{Optimizing with the operator's constraints, at industrial size}
\label{sec:opt-tool}
\label{sec:exp-scale}

In this section we set fidelity aside and use the model of
Section~\ref{sec:mip}, with the kernel calibrated on the analytical wake
of Section~\ref{sec:calibration}, as an optimization tool. The cell
is $\Delta = D/2$, as in the certificate experiments of~\cite{companion}
from which these campaigns descend; since every solver evaluates the same
implemented field, the comparison between solvers does not depend on how
faithful that field is, and the two heuristics of Section~\ref{sec:exp-big}
run unchanged at any cell.

\subsection{One model, three solvers, and the operator's constraints}
\label{sec:exp-scaling}
\label{sec:exp-heur}
\label{sec:exp-sidecon}

Three solvers consume the model in three forms. The commercial MILP
solver Gurobi~\cite{gurobi} solves as is the SOS2 formulation of
Section~\ref{sec:mip} on a ladder of square grids from $8\times8$ to
$20\times20$ cells, and closes only on some seeds of the smallest grid;
replacing its SOS2 condition by the upper concave envelope of the
breakpoints gives a relaxation, itself a MILP, that closes up to
$16\times16$\unskip, whose dual bound, even
at a time limit, bounds the exact objective of every feasible layout (\TN\
gives the proof) and certifies every layout of the ladder to within
8.75\%\ of the exact optimum. The second
solver, the \emph{vanilla local search} in what follows, is a greedy
construction followed by a one-exchange local search---single-turbine
relocation, first improvement, every candidate scored by the exact
$O(|\mathcal{C}|)$ evaluator of Section~\ref{sec:mip} with the exact power
curve---whose layouts are within 0.169\%\
of the MILP's on the $12\times12$ and $16\times16$ grids where the two
were compared (\TN). The third, the commercial heuristic solver
Hexaly~\cite{hexaly}, takes the binaries as the only decision variables
and obtains the field by propagation; it supplies no dual bound, hence it
serves as independent corroboration, the only check available at the
sizes of Section~\ref{sec:exp-big}. On the ladder, then, exact
optimization buys a declared gap and not a better layout. Beyond
$20\times20$ cells the ladder of certificates stops (the tests with the
noise constraint, \TN, reach $24\times24$): the model has
$O(|\mathcal{S}|\,|\mathcal{C}|)$ continuous variables, the certificate
stops closing on the ladder, and at the sizes of
Section~\ref{sec:exp-big} the tool is a heuristic whose only check is a
second, independent heuristic.

The constraints of Section~\ref{sec:related}---a cap on the machines of
each zone, exclusion areas of arbitrary shape, a site-dependent cost---are
each linear in the binaries and leave the propagation and the power
blocks untouched: a zone cap reads $\sum_{c \in \mathcal{A}_z} x_c \le
N_z$, an exclusion fixes $x_c = 0$ on its cells, and a site-dependent
cost enters the objective as $\sum_{c} C_c x_c$. A non-convex admissible
region---deep notches, annuli, a comb of corridors aligned with the
wind---leaves the exact solve and the local search within
0.005\%\ on the 5\
geometries we tried, and a free cardinality with a fixed cost per turbine
hurts the certificate and not the layout: the linear relaxation reverts to
a turbine on every free-stream site, the mechanism~\cite{companion} proves
for the uncosted model, and on one such instance the declared gap reaches
486\%\unskip.

What does change the picture is a constraint that couples all the binaries
at once. We take a noise immission limit $\sum_{c \in \mathcal{A}}
a_{rc}\, x_c \le L_r$ at receptors $r$ placed around the site, with
$a_{rc} = 1/\max(1, d(r,c)^2)$: a knapsack-like constraint that stands in
for a globally coupling condition, not a model of acoustic propagation.
On 12\ instances of a $20\times20$ grid,
with random receptors and severity and an equal budget per method, the
vanilla local search is the better of the two on average, but
2\ instances lose more than $5\%$, the
worst -11.3\%\unskip: there the local search
installs fewer turbines than the MILP, since a one-exchange neighborhood
cannot add one against a binding side constraint. The MILP fails for lack
of time and says so (largest declared gap
59.4\%\unskip), whereas the local
search fails silently, at a local optimum indistinguishable from success.

\subsection{Industrial sizes}
\label{sec:exp-big}

We finally give up the certificate and ask whether the formulation is of
use at industrial size, where what has to remain affordable is the
evaluation of a layout. On a $200\times200$\
grid---40\,000\ cells,
125~km$^2$ of site---evaluating a layout of
201\ turbines exactly takes
2.5~ms. On instances sized as farms,
Table~\ref{tab:biggrids} reports two heuristics, with
3\ independent seeds each and $600$~s per
run: the vanilla local search, extended with the simultaneous relocation
of two turbines and an install-or-remove move (\emph{vanilla LS++} for
short), and Hexaly. A layout of the largest instance is thus produced in
ten minutes on the laptop of Appendix~\ref{app:exp-conventions}; Valotta
Rodrigues et al.~\cite{valottarodrigues2024speeding} measure the
gradient-based route on PyWake at comparable sizes.

\begin{table}[htbp]
\centering\small
\caption{Two heuristics on the same model on large grids, \emph{vanilla LS++} (Section~\ref{sec:exp-big}) and Hexaly, at $600$~s per run and 3 seeds each, both scoring layouts with the same $O(|\mathcal{C}|)$ propagator. Value: expected power of the returned layout in kW, mean over the seeds; $n_t$: turbines installed; spread: relative distance between the best and the worst seed; gain: between the two means. Instances sized as farms (separation 5\,D, 1.6 turbines per km$^2$), one wind direction aligned with the grid at $10$\,m/s, $\Delta = D/2$. The lower panel charges a fixed cost per turbine and raises the cap to 2.5 times that of the upper panel; its objective is energy net of that cost, not comparable with the upper panel.}
\label{tab:biggrids}
\begin{tabular}{lrrrrrrrr}
\toprule
& & \multicolumn{3}{c}{vanilla LS++} & \multicolumn{3}{c}{Hexaly} & \\
\cmidrule(lr){3-5}\cmidrule(lr){6-8}
grid & $N_{\max}$ & value & spread & $n_t$ & value & spread & $n_t$ & gain \\
\midrule
\multicolumn{9}{l}{\emph{cardinality saturated}} \\[2pt]
100$\times$100 & 50 & 134\,949 & 0.1\% & 50 & 135\,051 & 0.0\% & 50 & +0.08\% \\
150$\times$150 & 113 & 302\,890 & 0.0\% & 113 & 303\,231 & 0.0\% & 113 & +0.11\% \\
200$\times$200 & 201 & 535\,856 & 0.0\% & 201 & 536\,906 & 0.1\% & 201 & +0.20\% \\
\midrule
\multicolumn{9}{l}{\emph{free cardinality, fixed cost per turbine}} \\[2pt]
100$\times$100 & 125 & 4\,616 & 2.8\% & 52 & 4\,777 & 1.3\% & 51 & +3.47\% \\
150$\times$150 & 282 & 8\,473 & 0.9\% & 101 & 8\,645 & 1.4\% & 96 & +2.03\% \\
200$\times$200 & 502 & 13\,331 & 1.6\% & 162 & 14\,286 & 0.9\% & 157 & +7.17\% \\
\bottomrule
\end{tabular}
\end{table} 

The upper panel corresponds to practice, where the tendering authority
fixes the installed capacity, hence the machine count, before any layout
is drawn. There the two searches are close: over the three grids Hexaly is
ahead of \emph{vanilla LS++} by between
0.1\%\ and
0.2\%\ on the mean of the seeds, and each is
stable across its own seeds to within
0.1\%\unskip: corroborating evidence, not
a proof of optimality. The lower panel charges a fixed cost per installed
turbine, equal to 96\%\ of the
expected power of one machine in free stream, and raises the cap to
2.5\ times that of the upper panel, so
that the count becomes a decision of the model. Here Hexaly is clearly
ahead, by up to 7.2\%\ on the mean of the
seeds, a gap that looks larger than it is in energy, since the objective
is energy net of the cost; Hexaly installs fewer machines on every grid
(columns $n_t$), so part of the margin is a decision on how many turbines
to install, not on where, which \emph{vanilla LS++} changes one machine
at a time. Under the twelve-sector rose of Horns Rev~1 (\TN), the seeds
return layouts of close but unequal value, spread over
0.27\%\ on the larger grid, which is why we
report spreads over seeds rather than a best-of. 
\section{Directional quadrature error and spatially localized surrogate bias}
\label{sec:sirens}
\label{sec:fidelity-not-opt}

In this section we measure on Lillgrund the two errors of a discretized
layout objective: the directional quadrature error, which draws a search
to the orientations between the directions of the rose (we nicknamed
these orientations the \emph{sirens} of the rose in an earlier version of
this work), and the spatially localized surrogate bias, which makes some
cells look freer than a RANS says they are (the \emph{mermaid
cells}). We compare three search pipelines, one on PyWake and two on the
EndoWake field at $\Delta = D$ and $\Delta = D/3$, at a RANS judge that
none of them calls. At randomly drawn directions the PyWake pipeline has
the largest mean gain, but five seeds resolve none of the paired margins.

\subsection{The test bench: three pipelines and a RANS judge}
\label{sec:rans-test}

We compare complete pipelines, as a planner would run them: in the sense
of Section~\ref{sec:exp}, the arms share the start, the budget, the power
curve and the calibration data, and differ in the wake model, in the wind
climate each integrates, in how each scans the neighborhood and in the
hardware, so the bench is a comparison between pipelines by design,
and not an ablation of the wake model alone. For each of five seeds,
every arm starts from the same certified local optimum of an
interference-matrix guide on the grid at $\Delta = D$, keeps its number of
machines, and runs the one-exchange local search of
Section~\ref{sec:exp-heur} for $1\,800$ seconds, restarting with a new
scan order at every local optimum and keeping the best layout found.
Every arm fits its own model to the same single-turbine RANS field by the
protocol of Section~\ref{sec:rans-calib}. Arm~A calls PyWake~2.6.20 with
\texttt{Bastankhah\_PorteAgel\_2014} at the fitted $k^{*} =
0.0337$, scores the annual energy under the
full Weibull distribution at 30\
directions, and moves each machine to the best of a sample of target
cells, evaluated in parallel on one CPU node. Arm~B runs the EndoWake
field at $\Delta = D$, with the one-cell correction of
Section~\ref{sec:fixpoint}, on a data-center GPU through the batch
evaluator of Section~\ref{sec:exp-batch}; it reads the twelve measured
sectors of the site, each integrated by five\
Gauss--Legendre nodes, and accepts the first improving move of a random
scan of the whole neighborhood. Arm~B$'$ is EndoWake at $\Delta = D/3$ on
the same sites, with the footprint of~\eqref{eq:mip-wake}.
Appendix~\ref{app:protocols} gives the instance, the sample sizes, the
stopping rules and the restart counts.

The judge is the RANS of Section~\ref{sec:rans}, run once on each final
layout at 8\,m/s over twelve directions drawn at random for each seed (the
\emph{random judge}): one uniformly in each thirty-degree interval between
consecutive sector centers (Appendix~\ref{app:rans-details}). The judged
farm power is the mean of the twelve powers, weighted by the frequency of
the sector that contains each direction, and the gain of a layout is taken
over the start judged at the same directions. The draw is independent of
every search, but 9\ of the
60\ drawn directions are nodes of the rose
of arm~A, and none of arm~B. No arm optimizes the functional the judge
measures, since each integrates its own rose, whereas the machine is
shared: all three arms use the tabulated power curve of the SWT-2.3-93,
and the judge its thrust and power tables. The bench thus also tests
whether a single-wake fit transfers to interacting machines, which is not
granted~\cite[Sect.~2.5]{thomas2023les}. A margin measured at 8\,m/s,
below rated, is a margin in farm power at that speed, and we neither
convert it into annual energy nor price it.

\paragraph{The site.}
\label{sec:exp-lillgrund}
Lillgrund is the offshore farm whose 48\
machines and row-wise SCADA deficits~\cite{sebastiani2021} are public, and
whose wakes have been measured by scanning lidar~\cite{sood2022lillgrund};
the wind enters as its 12-sector rose.
Public data do not give the available area, hence every arm shares the
convex hull of the built positions, less the two lattice sites the
designers left empty: the comparison is conditional on that area, and
for the same reason we do not compare any arm with the built layout.

\subsection{Directional quadrature error}
\label{sec:cones}

\begin{figure}[tbp]
\centering
\includegraphics[width=\linewidth]{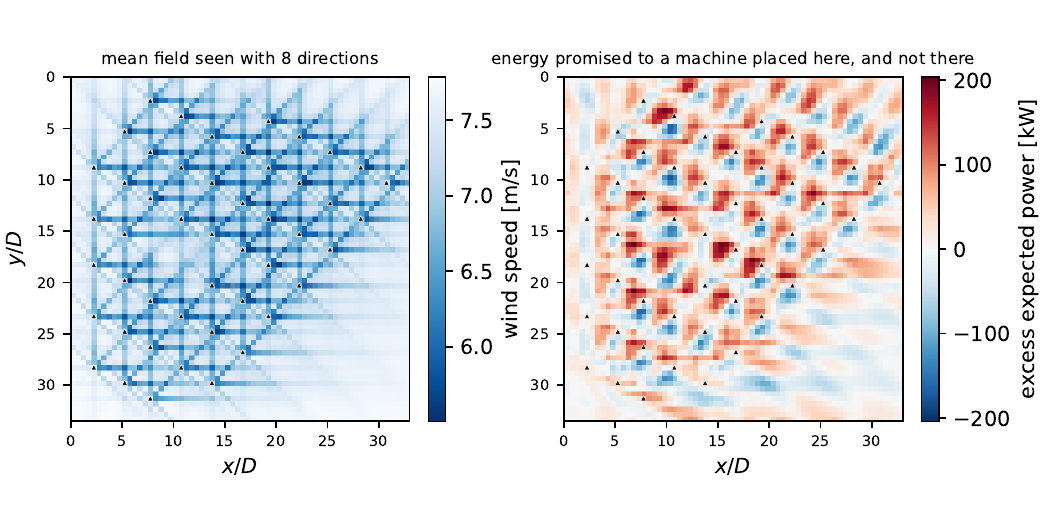}
\caption{Directional quadrature error on the site. Left: the
probability-weighted EndoWake field at $\Delta = D/2$ over the built
Lillgrund layout, under a rose of eight directions. Right: the expected
power this rose promises to a machine placed in each cell, minus what it
delivers once each sector is integrated by quadrature; red is power
promised and not there, blue is power delivered and not promised (averaged
over $1.5\,D$ for legibility).}
\label{fig:cones_field}
\end{figure}

The first defect is the discretization of the wind climate. A wake is
narrow next to a measured sector, so whether two turbines interfere
depends sharply on the angle between the wind and the line joining them.
A rose of $|\mathcal{S}|$ sectors offers the model $|\mathcal{S}|$
directions and no others, and a pair aligned between two of them is in
nobody's wake, for the model but not for the wind; a local search does
not average these orientations out but selects for them, since they are
where its objective is largest. The effect is known: with thirty-degree
bins, ``the optimizer smartly puts turbines in between the
cones''~\cite[Appendix~A.3]{fischetti2017thesis}, and it was reported
elsewhere~\cite{feng2015,quaeghebeur2021}, also against large-eddy
simulation~\cite{thomas2023les}. Figure~\ref{fig:cones_field} maps it on
our site under a deliberately coarse rose of
8\ directions: the power the rose promises to
a machine and does not deliver lies in the corridors between the wake
fans, and peaks at 15\%\ of its expected power.

A judge that samples the wind at the directions a search optimizes on
cannot measure the integration error of that search between its nodes.
The natural judge samples the six sector centers $0^\circ, 30^\circ,
\dots, 150^\circ$ with equal weights, where the central node of each
sector of arm~B, carrying 28\%\ of its
weight, falls exactly; arm~A optimizes over thirty directions twelve
degrees apart and shares at most three of them with this judge. At the six
directions halfway between the centers, $15^\circ, 45^\circ, \dots,
165^\circ$, with the same weights, arm~B gains
$-1.48\%$ where it gained
+4.06\%\ at the centers, and arm~A moves from
+1.13\%\ to +1.04\%\ (means over
the five seeds). We
call half the difference between the two judgments the \emph{phase
sensitivity} of an arm: it is 0.05\ points for
arm~A and 2.77\ for arm~B, and the margin of
arm~B over arm~A, 2.93\ points at the
centers, is $-2.52$ at the offsets. The two
judges share spacing and weights and differ only by a shift of fifteen
degrees, hence this change is an effect of phase alone. Phase sensitivity
needs twice the judging budget and no new search, and we propose it as a
first diagnostic for any tool that optimizes a discrete rose.

Arm~B takes the direction uniform within each sector and integrates it by
an $m$-node Gauss--Legendre rule, as in the quadrature estimates of the
energy yield of Padr\'on et al.~\cite{padron2019pce}, at the cost of $m$
replications of the field: on one pair of machines four diameters apart,
computed with PyWake alone, five nodes per sector cut the
orientation-averaged error of the bare twelve-sector rose from
1.79\ to 0.02\
percentage points (\TN). Quadrature does not remove the directional error
but moves it between the nodes: run with one node per sector, arm~B has a
phase sensitivity of 11.06\ points, which five
nodes shrink to 2.77\unskip.

\subsection{Spatially localized surrogate bias}
\label{sec:mermaid}
\label{sec:debias}

\begin{figure}[htbp]
\centering
\includegraphics[width=0.62\linewidth]{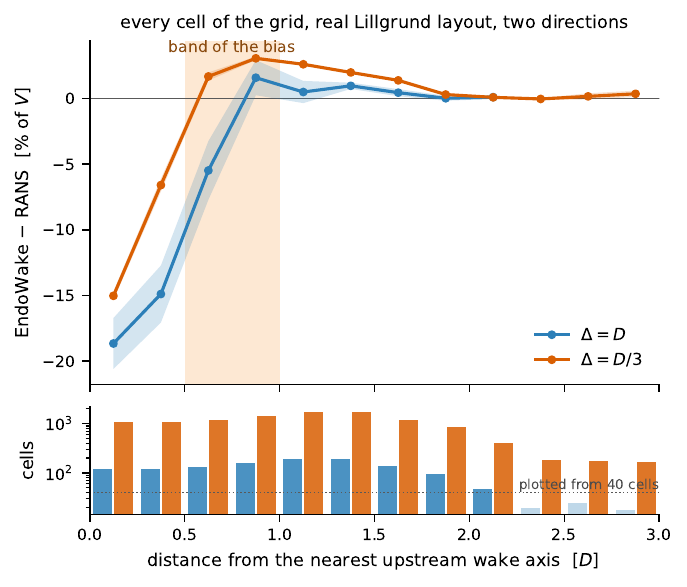}
\caption{The band of the spatially localized surrogate bias. Error of
the distilled field against the RANS, in percent of the free stream, on
every cell of the hub plane of the built farm at Lillgrund, binned by the
lateral distance of the cell from the axis of the nearest wake upwind, at
$\Delta = D$ (blue) and $D/3$ (orange), averaged over the directions
$222^\circ$ and $120^\circ$. Bins are $0.25\,D$ wide; the shading is one
standard error of the cell values, a descriptive dispersion since
neighboring cells are correlated; bins with fewer than 40 cells (light
bars in the lower panel) are not plotted. The shaded band, from $0.5$ to
$1\,D$, is the band of the text. Positive values mark overestimation, where
EndoWake promises more wind than the RANS delivers; negative values mark
underestimation.}
\label{fig:banda}
\end{figure}

\begin{figure}[tbp]
\centering
\includegraphics[width=0.9\linewidth]{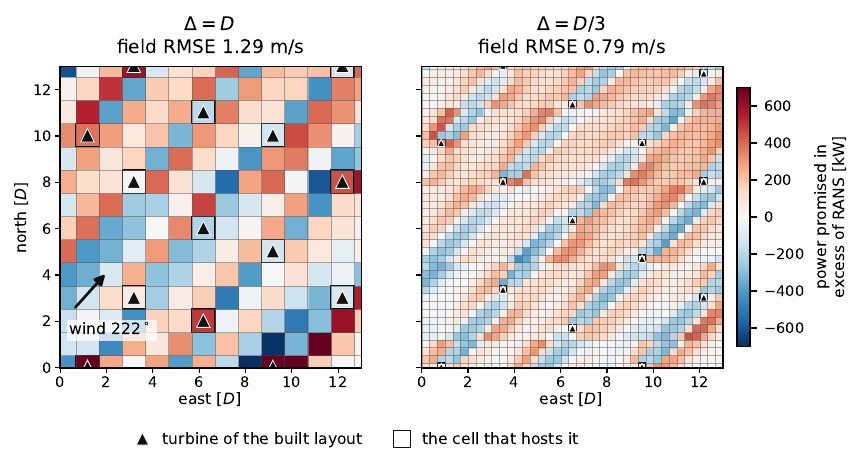}
\caption{Spatially localized surrogate bias on a window of the built
Lillgrund farm, wind from $222^\circ$ at 8\,m/s: for each cell, the power
EndoWake promises to a machine placed there minus what the RANS field of
the built layout, averaged over the cell, would give it, both read off
the existing fields with no machine added, at $\Delta = D$ (left) and
$D/3$ (right). Triangles are the built turbines; axes in diameters.
Red marks overestimation, where EndoWake promises more power than the
RANS delivers; blue marks underestimation.}
\label{fig:mermaid}
\end{figure}

The second defect is the grid, which does to space what the rose does to
direction. Our field, propagated at a fixed step, can reopen the flow just
outside the edge of a wake earlier than the RANS does, and be optimistic
relative to it in a band half a diameter to one diameter off the axis of
the machine upwind. Figure~\ref{fig:banda} measures this bias on every
cell of the hub plane of the built farm, independently of where any
optimizer put its machines. Both fields are pessimistic on the axis; the
error of the coarse field crosses zero 0.80\,D\ off
the axis and that of the fine field 0.55\,D\unskip,
i.e., the fine field becomes optimistic relative to the RANS closer to the
wake axis. At the machines of the bench the picture is the same: near a
wake edge ($0.5$ to $1.5$ diameters off the axis) the coarse field errs on
the pessimistic side ($-1.2\%$ of the free
stream) and the fine field promises +2.4\%\
that the RANS does not deliver, under winds along the grid axes and
oblique to them alike, whereas the error of the coarse field changes
sign between the two (\TN). Figure~\ref{fig:mermaid} shows the same bias
cell by cell under one wind: the coarse field scatters its error over the
grid, whereas the fine field, closer to the RANS on average, resolves the
corridors between the wakes and lines their edges with optimistic cells.

Unlike the directional error, this bias is present under a single wind,
and the rose only sums it. As a hypothesis, not a measured mechanism, we
attribute it to the turbulent mixing that widens the wakes and fills the
corridors between them, which a linear, local propagation calibrated on
one isolated wake cannot represent. A planner can locate this bias on a
new site with a few RANS directions on one layout, binned as in
Figure~\ref{fig:banda}. Whether a search exploits it is another matter:
averaged over the directions of the rose,
26.0\%\ of the machines of the start sit in the
band, and the optimized layouts change that share by
$-0.5$ points in arm~B,
$-0.6$ in arm~B$'$ and
$0.1$ in arm~A. On this site, the band is
where the field is wrong, not where its search goes.

\subsection{Results at the random judge}
\label{sec:sirens-judge}

Table~\ref{tab:rans_seeds} reports, seed by seed, the gain of each arm
over the start at the random judge of Section~\ref{sec:rans-test}. Paired
seed by seed, the margin of arm~B over arm~A averages
$-0.78$ points (standard error
0.57\unskip, 95\% $t$-interval
$[-2.36,
0.80]$), and that of arm~B$'$
$-0.80$ (standard error
0.94\unskip, interval
$[-3.40,
1.80]$). Both intervals contain zero
and are wide: five paired starts detect no difference, and permit
materially different outcomes, so the order of the means is descriptive
and establishes neither a superiority of PyWake nor that refining the
field worsens the layouts in general.

\begin{table}[htbp]
\centering\small
\caption{Gain over the start, in percent of farm power, at the random RANS
judge (8\,m/s, twelve directions drawn for each seed, SWT-2.3-93), for
each of the five starts and on average. A: PyWake; B: EndoWake at $\Delta
= D$; B$'$: EndoWake at $\Delta = D/3$ on the same sites. Each row judges
the three pipelines on the same directions; the spread over rows is paired
variation across starts, not uncertainty of the RANS closure.}
\label{tab:rans_seeds}
\begin{tabular}{@{}lrrr@{}}
\toprule
Start & A: PyWake & B: EndoWake, $D$ & B$'$: EndoWake, $D/3$ \\
\midrule
0 & $+3.87$ & $+1.39$ & $-0.29$ \\
1 & $+1.36$ & $+1.13$ & $+1.86$ \\
2 & $+0.90$ & $+0.48$ & $-0.45$ \\
3 & $+0.46$ & $-1.12$ & $+1.62$ \\
4 & $+0.21$ & $+1.00$ & $+0.08$ \\
\midrule
Mean & $+1.36$ & $+0.58$ & $+0.56$ \\
\bottomrule
\end{tabular}
\end{table}

Three sources of error meet here. The margin of arm~B over arm~A changes
sign between the sector centers and the offsets (Section~\ref{sec:cones}),
is nearly zero at their union (+1.08\%\ for arm~A
against +1.25\%\ for arm~B) and
$-0.78$ at the random judge, which covers the
whole circle with the weights of the rose: the ranking is sensitive to the
evaluation. The layouts of arm~B$'$ do not populate the band of
Section~\ref{sec:mermaid} more than the start, so the spatial bias does
not by itself explain a mean below that of arm~A. And the twelve
directions of the judge are a sample, redrawn for each seed, on a
reference with its own closure and errors (Section~\ref{sec:rans}). The
bench spends the throughput of the distillate on more descents of the
same search from one start; what that throughput is worth to a search
built for it, Section~\ref{sec:proxy} measures.
That search also gives the ablation the bench cannot. Two of its
arms on the unseen seeds, PyWake alone and EndoWake alone, run the same
search from the same start, with the same neighborhood, acceptance rule
and stop rule, and differ in the evaluator only: PyWake on thirty
processes, or the distillate on one GPU, accepts the moves, and the
hardware differs with the evaluator, as in Table~\ref{tab:env}. We
judged the final layouts of the first 5%
\
unseen seeds (5 to 9%
) at the random judge
of Section~\ref{sec:rans-test}, 12%
\
directions drawn for each seed, 240%
\ RANS
runs in all, the starts included (\TN). Seed by seed, with the weights
of the rose, the gain over the common start at the judge is
3.64\%%
\ for PyWake alone and
3.30\%%
\ for EndoWake alone, a paired
difference of $-$0.34%
\ points (standard
error 0.34%
\unskip, 95\% $t$-interval
$[$-$1.29%
,
$+$0.61%
]$), where PyWake reads
$-$0.21%
\ points between the same layouts.
The third arm, screen-and-confirm, judged on the same directions, gains
3.79\%%
\unskip,
$+$0.16%
\ points from PyWake alone
($[$-$0.96%
,
$+$1.27%
]$), where PyWake reads
$+$0.03%
\ points. Both intervals contain
zero and still permit materially different outcomes: at the reference,
the layouts the distillate chooses on its own, and those it screens for
PyWake, are not told apart from those PyWake chooses, within what five
paired starts resolve, and PyWake reads the same small differences
between them. 
\section{What the throughput buys: EndoWake as a proxy of PyWake}
\label{sec:proxy}

In this section we use EndoWake, fitted on PyWake, to speed up a local
search whose objective is PyWake itself, so that no RANS is needed:
EndoWake ranks a sample of the neighborhood on a GPU, on a reduced
description of the wind, and PyWake confirms the head of the ranking on
the full objective. On twenty held-out starts on Lillgrund, this
screen-and-confirm search and the search on PyWake alone reach similar
mean PyWake-valued gains, and the former stops about an order of
magnitude sooner at one wind speed and on the annual energy, less when
the proof of local optimality is charged to the timed search
(Table~\ref{tab:proxy}).

\paragraph{The instance.} The farm is Lillgrund as in
Section~\ref{sec:exp-lillgrund}, with sites at $\Delta = D/3$
(4\,728\ allocatable cells) and a rose of $360$
directions one degree apart at 8\,m/s, each carrying a thirtieth of the
frequency of its sector. We fit EndoWake by the recipe of
Section~\ref{sec:calibration} on the single-wake field of the Gaussian
model of PyWake at $k = 0.0345$, with a residual of
6.8\%\unskip, and let it read the tabulated power
curve of the SWT-2.3-93, so that the two models share the machine.
This $k$ is a fixed parameter of the objective in every search of this
section: it is the value fitted on the RANS before the correction of
Section~\ref{sec:rans}, and it differs from the $k^{*}$ of
Table~\ref{tab:rans_calib}. The comparison between searches does not rest
on it. Evaluated at $k^{*}$, the final layouts on the unseen seeds gain
0.05\ points more at 8\,m/s and
0.04\ on the annual energy, in every arm, and
the paired differences from PyWake alone move by at most
0.002\ points. The
starts are twenty-five layouts of 44--48\
machines: the five of Section~\ref{sec:rans-test} (the development seeds
$0$ to $4$), on which we designed the search, and twenty new starts on the
same site (the unseen seeds $5$ to $24$), built the same way and run only
after the configuration and the policy of the list had been frozen. A
move relocates either one machine, to any allocatable cell that keeps the
minimum separation, or two machines at most $5\,D$ apart, each to a cell
within $2\,D$ of its position, in the spirit of Cazzaro and
Pisinger~\cite{cazzaro2022vns}. Counting all the generated moves, the
neighborhood holds 37\ to
79\ thousand moves at the starts, and up to
145\ thousand along a search at 8\,m/s
(153\ thousand over the annual-energy
runs below). The neighborhood is the same in every arm; PyWake alone
enumerates it again after every accepted move, screen-and-confirm after
every round with an accepted move (Figure~\ref{alg:screen}), and every
arm runs from the start to its first stop, with no restart
(Appendix~\ref{app:sm-proxy}).

\paragraph{The arms.} (i)~\emph{PyWake alone} scans the neighborhood in
random order, in blocks evaluated in parallel, and accepts the first
improving move, as in Section~\ref{sec:exp-heur}. Best improvement, as in
Cazzaro and Pisinger~\cite{cazzaro2022vns}, is out of reach: their
interference matrix prices a move in a number of sums linear in the
number of machines, whereas one full PyWake scan of our neighborhood costs
about 43~s, and PyWake alone accepts
155\ moves per run. (ii)~\emph{EndoWake
alone} runs the same search on the distillate, and PyWake evaluates its
final layout. (iii)~\emph{Screen-and-confirm} is a surrogate-ranked
candidate selection with PyWake acceptance: EndoWake orders the moves,
and PyWake accepts only those that improve its own value
(Figure~\ref{alg:screen}).

\paragraph{Screen-and-confirm.} Its design follows from the time scales
of the two models. EndoWake evaluates
23\ times as many layouts per second
as PyWake on thirty processes, but ranking one GPU batch of $8\,192$ moves
takes 0.26~s, the time of
358\ PyWake evaluations, whereas PyWake
alone finds each of its first ten moves in
31\ evaluations on average: a ranking
before every move would cost about
11\ times the move it serves.
We therefore reuse each ranking over a round of moves, since a move
displaces one or two machines out of 44--48\unskip,
and start each round from a random sample of one batch, doubled only
after a round without accepted moves, since improving moves are plentiful
early in a search and rare close to the local optimum. A \emph{round} is one pass of
PyWake down the first $L$ (say) moves of the ranking, each evaluated on
the current incumbent, hence every acceptance is an exact improvement of
the PyWake value, as in the search on PyWake alone; several moves can be
accepted in one round, provided they touch different machines, and only
the order of the proposals ages along it. After a round with an accepted
move the neighborhood is enumerated again, the sample is redrawn and
ranked afresh, and the set of machines moved in the round is cleared; the
record of rejected moves is kept, but keyed to the incumbent at which
each rejection was observed. EndoWake ranks in single precision, as it only orders
the moves. The scheme is a candidate list in the sense of tabu
search~\cite{glover1997tabu}, and a multifidelity scheme in which, as in
Jasa et al.~\cite{jasa2022multifidelity}, a low-fidelity wake model guides
and a higher-fidelity one corrects.

\paragraph{The list and the stop.} The length $L$ of the list is not
tuned. At its start, inside the clock, the search measures the rate $r_P$
of PyWake on $512$ random neighbors of the start and the rate $r_E$ of
EndoWake on one batch, and it sets $L = \lceil q\,r_P / r_E \rceil$ for a
sample of $q$ moves, so that PyWake spends on confirming about the time
EndoWake spends on ranking. Over whole runs on thirty processes,
confirming took 1.3\ to
1.9\ times as long as ranking.
The search stops after a round without accepted moves on the whole
neighborhood, at no extra cost. PyWake has then seen only the head of the
last ranking, so the stop is not a proof of local optimality; the
\emph{certified} variant adds, inside the clock, a scan of the whole
final neighborhood by PyWake, resumes the search if the scan finds an
improving move, and stops only when none is left, so that its final
layout is a local optimum of PyWake with the same proof as the search on
PyWake alone. For the uncertified runs we check the final layouts
afterwards, off the clock.

\paragraph{The screen.} EndoWake need not read the wind of the objective,
since it only orders the moves and PyWake confirms each of them on the
full objective. When the objective integrates several wind speeds,
EndoWake ranks on a reduced description of the wind: a few speed nodes
below the rated speed of the machines, each carrying, direction by
direction, the probability of the speed bins around it, and no node above,
where the machines run close to their rated power and their wakes cost
little. The reduction concerns the ranking only: PyWake confirms
every move on all the speed bins of the objective, so the wake losses at
high speeds are charged in full to every accepted move, and the
annual-energy runs below compare this screen with one that reads all the
speeds. At a single speed, as in the instance above, the screen is the
wind of the objective itself, and the method reduces to ranking on the
rose of the objective.

\begin{figure}[htbp]
\centering
\fbox{\begin{minipage}{0.94\linewidth}\small
\emph{Input:} a start layout $x$ with its neighborhood $\mathcal{N}(x)$, the PyWake
value $F$ on the objective, the EndoWake score $E$ on the screen, and one GPU
batch $q_0 = 8\,192$.
\begin{enumerate}
\setlength{\itemsep}{1pt}
\item Inside the clock, measure the evaluation rates $r_P$ of $F$ and $r_E$
of $E$ on random neighbors of $x$.
\item Enumerate $\mathcal{N}(x)$, draw a random order of it, and set
$q = \min\{q_0, |\mathcal{N}(x)|\}$.
\item Score by $E$ the moves among the first $q$ of the order that have no
score yet, rank these $q$ moves by score whatever its sign, and set
$L = \lceil q\, r_P / r_E \rceil$.
\item Go down the first $L$ moves of the ranking. Skip a move that is no
longer feasible, touches a machine moved in this round, or was already
rejected at this incumbent (a move rejected at an earlier incumbent
is evaluated again); evaluate every other move by $F$ on the current
$x$, and accept it if it improves $F(x)$.
\item If a move was accepted, go to Step~2. If $q < |\mathcal{N}(x)|$, set
$q = \min\{2q, |\mathcal{N}(x)|\}$ and go to Step~3.
\item Stop: return $x$. Certified variant: scan $\mathcal{N}(x)$ by
$F$; if an improving move is found, accept it and go to Step~2, else
return $x$.
\end{enumerate}
\end{minipage}}
\caption{Screen-and-confirm. PyWake evaluates the moves in parallel
blocks, one move per process (Appendix~\ref{app:sm-proxy}). A
rejection is tied to the incumbent at which it was observed: a move is
never skipped only because an earlier incumbent rejected it.}
\label{alg:screen}
\end{figure}

\paragraph{Hardware and timing.} PyWake runs, in every arm, one
simulation per process on 30\ worker
processes of one node; EndoWake runs on one GPU of the same node, in
batches of $8\,192$ layouts carried through all the directions of the rose
in one call (Section~\ref{sec:exp-batch}); processor, GPU and software
versions are in Appendix~\ref{app:exp-conventions}. Screen-and-confirm
thus adds one GPU to the processes of PyWake alone
(Section~\ref{sec:exp-batch} compares the evaluators at equal hardware).
All times are warmed search
wall times: the warm-up of the two models, the just-in-time compilation,
the start of the process pool and the fit of the distillate are off the
clock, whereas the rate measurement is on it. It measures $r_P
=1\,369$ and $r_E =
31\,393$ layout evaluations per second on
average, hence a list of 4.4\% of the
sample. To see the list follow the machine, we repeat the unseen seeds
with PyWake on eight processes, the CPU share of one GPU on our cluster.

\begin{table}[htbp]
\centering\small
\caption{EndoWake as a proxy of PyWake. Lillgrund, sites at $\Delta = D/3$; the
objective is the farm power under $360$ directions at 8\,m/s, or the expected
power under the annual wind of the site ($360$ directions by 22 speed bins);
PyWake (\texttt{Bastankhah\_PorteAgel\_2014} at $k = 0.0345$, SWT-2.3-93) is the
objective and the judge. Every arm runs the same neighborhood from the same start
to its stop, with no restart, on twenty unseen seeds. Screen-and-confirm ranks on
the wind of the objective at 8\,m/s and on four speeds below rated for the annual
energy. With the certificate: the same search followed, inside the clock, by a
scan of the whole final neighborhood by PyWake, resumed if the scan finds an
improving move. Gain: power, or expected power, of the final layout over the
start, by PyWake, in percent. delta: paired mean difference from PyWake alone, in
percentage points (standard error over the seeds in parentheses). $t$: warmed
wall-clock seconds, rate measurement included; cert.: seconds of $t$ spent by the
final scans. Speedup: over PyWake alone, geometric mean. Local
optima: final layouts from which PyWake alone finds no improving move. The
certified runs are separate runs, and at 8\,m/s their final gains differ from
those of the uncertified ones by up to 0.40 points. TOPFARM: SLSQP on the
gradients of PyWake from the same starts, under the separation, the hull and the
two forbidden disks, on one process, with no neighborhood, in the plain
configuration or continued on the boundary with the best feasible recorded
point.
Length of the list: screen-and-confirm rerun with the list multiplied by the
factor, or ranked afresh after every accepted move. Horns Rev~1: the second site of Section~\ref{sec:proxy}.}
\label{tab:proxy}
\begin{tabular}{lrrrrrc}
\toprule
 & gain & delta & $t$ & cert. & speedup & local optima \\
\midrule
\multicolumn{7}{l}{\emph{8\,m/s, unseen seeds, PyWake on 30 processes}} \\
\quad PyWake alone & $+4.33$ &  & 294 &  &  & 20/20 \\
\quad screen-and-confirm & $+4.28$ & $-0.05$\,(0.07) & 38 &  & 7.9$\times$ & 19/20 \\
\quad \quad with the certificate & $+4.30$ & $-0.04$\,(0.07) & 98 & 58 & 3.0$\times$ & 20/20 \\
\quad TOPFARM, continued SLSQP & $4.30$ & $-0.03$\,(0.26) & 54 &  &  &  \\
\addlinespace
\multicolumn{7}{l}{\emph{8\,m/s, unseen seeds, PyWake on 8 processes}} \\
\quad PyWake alone & $+4.33$ &  & 945 &  &  & 20/20 \\
\quad screen-and-confirm & $+4.32$ & $-0.01$\,(0.06) & 39 &  & 24.5$\times$ & 17/20 \\
\quad \quad with the certificate & $+4.28$ & $-0.05$\,(0.07) & 240 & 200 & 4.1$\times$ & 20/20 \\
\addlinespace
\multicolumn{7}{l}{\emph{8\,m/s, length of the list, 20 unseen seeds, PyWake on 30 processes}} \\
\quad list $\times 0.5$ & $4.29$ & $-0.04$\,(0.05) & 26 &  &  & 17/20 \\
\quad list $\times 1$ (rerun) & $4.28$ & $-0.05$\,(0.06) & 37 &  &  & 19/20 \\
\quad list $\times 2$ & $4.35$ & $+0.02$\,(0.07) & 57 &  &  & 20/20 \\
\quad list $\times 4$ & $4.31$ & $-0.03$\,(0.07) & 98 &  &  & 20/20 \\
\quad one move per ranking & $4.30$ & $-0.03$\,(0.06) & 45 &  &  & 19/20 \\
\addlinespace
\multicolumn{7}{l}{\emph{annual energy, unseen seeds, PyWake on 30 processes}} \\
\quad PyWake alone & $+2.94$ &  & 2\,421 &  &  & 20/20 \\
\quad screen-and-confirm, four speeds & $+2.90$ & $-0.04$\,(0.04) & 161 &  & 14.9$\times$ & 18/20 \\
\quad \quad screen on all 22 speeds & $+2.90$ & $-0.04$\,(0.04) & 854 &  & 2.9$\times$ & 20/20 \\
\quad \quad with the certificate, all speeds & $+2.90$ & $-0.04$\,(0.04) & 1\,183 &  & 2.0$\times$ & 20/20 \\
\quad \quad with the certificate, four speeds & $+2.89$ & $-0.05$\,(0.04) & 496 & 333 & 4.7$\times$ & 20/20 \\
\quad TOPFARM, plain SLSQP & $+2.54$ & $-0.40$\,(0.10) & 391 &  &  &  \\
\quad TOPFARM, continued SLSQP & $+2.84$ & $-0.10$\,(0.18) & 315 &  &  &  \\
\addlinespace
\multicolumn{7}{l}{\emph{Horns Rev~1, 8\,m/s, sites at $\Delta = D$, 10 seeds, PyWake on 30 processes}} \\
\quad PyWake alone & $+9.20$ &  & 2159 &  &  & 10/10 \\
\quad screen-and-confirm & $+8.91$ & $-0.29$\,(0.04) & 65 &  & 33.2$\times$ & 4/10 \\
\quad \quad with the certificate & $+8.96$ & $-0.24$\,(0.04) & 379 & -- & 5.7$\times$ & 10/10 \\
\bottomrule
\end{tabular}
\end{table} 

Table~\ref{tab:proxy} reports the outcome. At 8\,m/s, on the unseen seeds,
PyWake alone gains +4.33\%\ in
294~s on average, with
391\ thousand PyWake evaluations per
run. Screen-and-confirm gains +4.28\%\ in
38~s, the rate measurement
(0.6~s) included, with
20\,277\ PyWake and
422\ thousand EndoWake
evaluations. A run accepts 90\
moves on average, 10\ of them on
its first ranking and 4.2\ per
ranking that yields any, against one per scan for PyWake alone. The mean
paired difference from PyWake alone is
$-$0.05\ points, with a 95\% $t$-interval
of [$-$0.19,
$+$0.10] over the twenty seeds, and the
speedup ranges from 4.1\ to
14.4\ over the seeds. The two
searches end on different local optima, hence differences of either sign
from seed to seed, from $-$0.78\ to
$+$0.40\ points. In
19/20\ seeds the final layout is also
a local optimum of PyWake: PyWake alone, started from it off the clock,
finds no improving move; in the remaining seed it finds
23\ more moves, worth
0.17\ points. The certified
runs, which carry that proof inside the clock, end in
98~s on average, of which
58~s go to the final scans: a speedup
of 3.0\ (from
1.2\ to
6.1\ over the seeds) at a paired
difference of $-$0.04\ points
([$-$0.18,
$+$0.11]), with every final layout a
local optimum of PyWake.

\begin{figure}[htbp]
\centering
\includegraphics[width=0.85\linewidth]{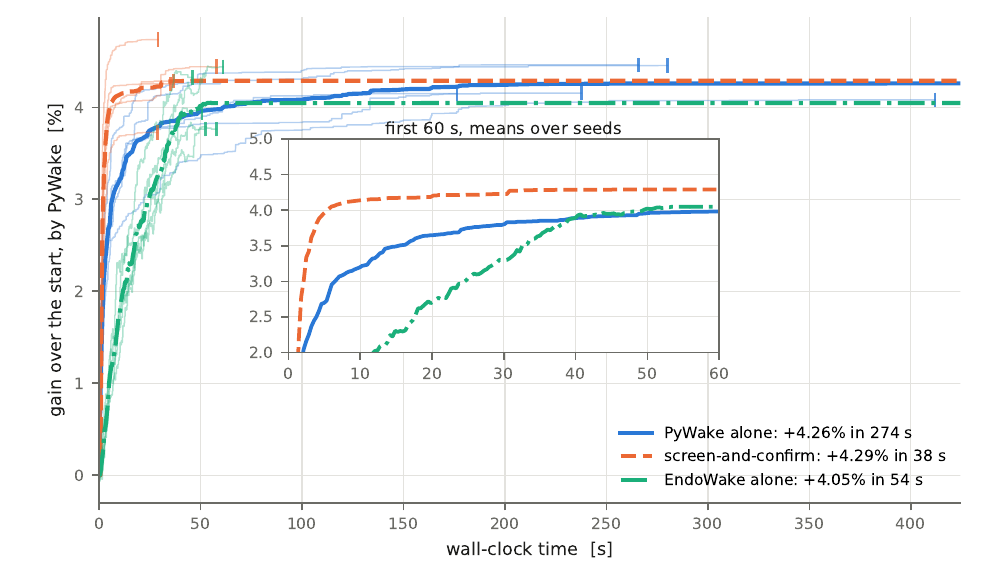}
\caption{Gain against time, unseen seeds, PyWake on thirty processes:
gain of the incumbent over the start, by PyWake, against wall-clock
seconds. Thick lines are means over seeds $5$ to $9$, the only ones run
also with EndoWake alone, whose curve is the PyWake value of its
trajectory computed off the clock; thin lines are the seeds, with a tick
where each search stops; the inset shows the first minute.
Screen-and-confirm accepts its first move
\protect0.9~s into the run, rate
measurement included, and leads from
\protect1.4~s on.}
\label{fig:proxy-time}
\end{figure}

Figure~\ref{fig:proxy-time} shows where the time goes. Early in a search,
PyWake alone finds a move in a few tens of evaluations, and the sample of
one batch keeps the proxy from starting behind. The gap widens as
improving moves become rare: PyWake alone then spends thousands of
evaluations per move, and a full scan to declare the local optimum, while
one ranking still serves a whole round. EndoWake alone (dash-dotted curve)
ends below both on average, within two standard errors of PyWake alone
(Table~\ref{tab:proxy-roses}). The list follows the machine with no
intervention: with PyWake on eight processes, $r_P$ drops to
385\ evaluations per second while $r_E$
does not move, and the list shrinks from
4.4\% to
1.2\% of the sample. Screen-and-confirm
then takes 39~s, about as on thirty
processes, against 945~s for PyWake alone:
a speedup of 24.5\unskip, at a mean
difference of $-$0.01\ points
([$-$0.13,
$+$0.12]).

The development seeds say where the speedup comes from. Ranking the whole
neighborhood in double precision after every accepted move keeps the
quality (0.00\ points) but gives no speedup
(a ratio of 1.0\ to PyWake alone), since
the rankings take 98\%\ of the time. One
ranking per round, with a list of $3\,000$ moves tuned on these seeds,
brings the ratio to 4.5\unskip; single
precision to 6.9\unskip; the sample to
9.3\unskip, for
30\ CPU processes against one GPU.

\begin{table}[htbp]
\centering\small
\setlength{\tabcolsep}{4pt}
\caption{Directional quadrature error, with and without the confirmation. EndoWake reads a rose of $360$, $120$, $36$ or $12$ directions; PyWake, the objective and the judge, always reads $360$. Lillgrund at $\Delta = D/3$, unseen seeds $5$--$9$, PyWake on thirty processes, EndoWake on one GPU. Gain: power of the final layout over the start, by PyWake, in percent. delta: paired mean difference from PyWake alone (+4.26\% in 274~s), in percentage points, with its standard error in parentheses. $t$: wall-clock seconds. Speedup: over PyWake alone, geometric mean. Left: improving moves that PyWake alone still finds from the final layout, mean over the seeds. List: length of the list in percent of the sample.}
\label{tab:proxy-roses}
\begin{tabular}{lrrrrrrrrrr}
\toprule
 & \multicolumn{4}{c}{EndoWake alone} & \multicolumn{6}{c}{screen-and-confirm} \\
\cmidrule(lr){2-5}\cmidrule(lr){6-11}
directions & gain & delta & $t$ & speedup & gain & delta & $t$ & speedup & left & list \\
\midrule
360 & $+4.05$ & $-0.21$\,(0.14) & 54 & 5.0$\times$ & $+4.29$ & $+0.03$\,(0.13) & 38 & 7.3$\times$ & 5 & 4.6 \\
120 & $+3.70$ & $-0.56$\,(0.19) & 19 & 14.3$\times$ & $+4.21$ & $-0.05$\,(0.13) & 17 & 15.7$\times$ & 2 & 1.5 \\
36 & $+2.91$ & $-1.35$\,(0.07) & 7 & 38.3$\times$ & $+3.98$ & $-0.28$\,(0.23) & 8 & 32.9$\times$ & 42 & 0.4 \\
12 & $+0.41$ & $-3.85$\,(0.34) & 4 & 64.3$\times$ & $+1.38$ & $-2.88$\,(0.56) & 2 & 149.5$\times$ & 120 & 0.2 \\
\bottomrule
\end{tabular}
\end{table} 

\paragraph{Directional quadrature error, with and without the confirmation.}
Table~\ref{tab:proxy-roses} repeats the experiment on the unseen seeds
$5$ to $9$ with EndoWake reading a sparser rose, down to the twelve sector
centers, while PyWake keeps its $360$ directions. Without the
confirmation, the paired difference from PyWake alone is
$-$0.21\ points with the full rose,
$-$0.56\ with $120$ directions,
$-$1.35\ with $36$ and
$-$3.85\ with twelve, and only the first is
within two standard errors of zero. At twelve, the directional quadrature
error of Section~\ref{sec:cones} takes almost the whole gain: judged by
PyWake, those layouts are worth +0.41\%\ over
the start. Under screen-and-confirm, where PyWake rejects the moves that
only the sparse rose likes, the difference is
$+$0.03\unskip,
$-$0.05\ and
$-$0.28\ points with $360$, $120$ and $36$
directions, and $-$2.88\ with twelve: the
confirmation cures most of the loss down to $36$ directions, but not at
twelve. A sparser rose also shortens the list, since EndoWake ranks
faster (from 4.6\% of the sample to
0.2\% at twelve directions): a poorer
ranking meets a shorter list, the stop comes early, and PyWake alone
still finds 120\ improving moves from
the final layouts at twelve directions. The rule of the list knows
nothing of the quality of the ranking, hence EndoWake reads all the
directions of the objective in Table~\ref{tab:proxy}, also where the
screen reduces the speeds.
\paragraph{What the ranking sees.} The list is short against the
neighborhood, so the ranking must put improving moves at its head. On the
whole neighborhood of the starts of seeds $5$ to $9$, scored by PyWake off
the clock, 19\%%
\ of the
59%
\ thousand moves improve the PyWake value.
EndoWake ranks here the whole neighborhood, a stricter test than the
sample of one batch that the search ranks: of the first $720$ moves of
its ranking in single precision,
100.0%
\% do, of the first $1\,440$
99.5%
\%, and of the first $3\,000$
96.5%
\%, and the best move of the
neighborhood by PyWake sits at rank 3%
\
of the ranking on average. This is a diagnostic at the initial incumbent,
where improving moves are plentiful and the head of the ranking is thus
almost entirely made of them; whether the length of the list matters
along a search is measured next.
\paragraph{The length of the list.} Table~\ref{tab:proxy} also
reports the search rerun on the 20%
\
unseen seeds with the list of the rule, halved, doubled and quadrupled,
and with a fresh ranking after every accepted move, each run continued off
the clock by PyWake alone down to a local optimum of PyWake. The gains lie
between 4.28\%%
\ and
4.35\%%
\unskip, within the paired noise of
PyWake alone on the same seeds, whereas the search time grows with the
list, from 26%
~s at a mean list of
391%
\ moves to
98%
~s at 3300%
.
The continuation accepts 0.7%
\ moves on
average with the list halved (17%
\ of
20%
\ runs had stopped at a local optimum of
PyWake), 1.1%
\ with the list of the rule
(19%
\ runs), none with the list doubled or
quadrupled, and 1.0%
\ with a fresh ranking
after every move (19%
\ runs). Averaged over
the seeds, these moves are worth at most
0.01%
\ points of gain. These
20%
\ starts resolve no difference in quality
between the lists; a longer list
costs time, and the rule sits at the fast end of the range.

\paragraph{Annual energy.} We repeated the comparison on the objective a
planner optimizes, the expected farm power under the annual wind of the
site, which is proportional to the annual energy. PyWake integrates the
farm power over 360\ directions one degree
apart and 22\ speed bins of 1\,m/s between
3\ and 25\,m/s,
each weighted by the Weibull distribution of its sector at Lillgrund (the
site shipped with PyWake, fitted on seven months of measurements, on which
the annual energy is conditional), and both models read the tabulated
power and thrust curves of the SWT-2.3-93. We fitted EndoWake on PyWake at
the center of every bin, with residuals from
6.7\%\ to 13.1\%\unskip.
The screen has four\ nodes, at
5.5, 7.5, 9.5 and 11.5\,m/s, and none from
13\,m/s up. The starts are
twenty\ layouts of
43--48\ machines (seeds
10\ to
29\unskip: those up to $24$ are the unseen
starts above, the last five are new), and the protocol, the screen
included, was fixed before any search on the annual energy
(Appendix~\ref{app:sm-proxy}). On the full objective, PyWake evaluates
221\ layouts per second on thirty processes;
EndoWake evaluates 8\,038\ on the screen,
36\ times as many, against
1\,481\ (7\
times) when it reads all the 22\ speeds,
and the list shrinks from 14.8\% to
2.8\% of the sample.

PyWake alone gains +2.94\%\ in
2\,421~s on average (Table~\ref{tab:proxy}),
from an expected farm power of
34.2\ to
35.2~MW, i.e., from
300\ to
308~GWh of annual energy.
Screen-and-confirm gains +2.90\%\ in
161~s, a speedup of
14.9\ in geometric mean (95\% interval
11.9\ to 18.8),
at a paired difference of $-$0.04\ points
([$-$0.13,
$+$0.06]; from
$-$0.42\ to $+$0.24\
seed by seed). In 18/20\ seeds its final
layout is a local optimum of PyWake; in each of the other
two\unskip, a scan of the whole final
neighborhood, which leaves the layout as it is, finds
one\ improving move, worth at most
0.005\ points. Ranking on all the
22\ speeds gives a similar difference
($-$0.04\ points,
[$-$0.12,
$+$0.04]) at a speedup of
2.9\ only; the coarser screen gains about as
much on average, at a paired difference of $+$0.01\
points from the full one ([$-$0.06,
$+$0.07]), and the confirmation keeps every
accepted move exact. The certified runs, with the same reduced
screen and the final scans inside the clock, end in
496%
~s against
2\,421~s for PyWake alone, a speedup of
4.7%
\ (interval
4.0%
\ to
5.6%
) at a paired difference of
$-$0.05%
\ points
([$-$0.13%
,
$+$0.04%
]); the scans take
333%
~s of that time, and an exhaustive
scan by PyWake off the clock confirms every final layout as a local
optimum (\TN\ lists the runs). Ranking on all the speeds, the certified
runs are 2.0\ times faster only
(Table~\ref{tab:proxy}). The advantage does not come only from the stop: all
twenty\ starts reach a gain of $2\%$ in both
arms, screen-and-confirm 3.09\ times sooner in
geometric mean (95\% interval 2.25 to
4.24), whereas at $0.5\%$ the fixed cost of
the screen makes PyWake alone faster (0.55).
Re-evaluating the start and the final layout of every run with the step
of the quadrature halved in direction, in speed or in both, or with the
rose rotated by half a degree, moves its gain by up to
0.049\ points, but the mean paired difference
of screen-and-confirm from PyWake alone by at most
0.003\ points.

PyWake defines the objective of this experiment and evaluates every final
layout: the experiment measures how efficiently the screening optimizes
that objective, discretized as stated, not the physical accuracy of
either model; every error of EndoWake relative to PyWake is charged in
full, and no error of PyWake is seen.
\paragraph{A gradient-based baseline.} TOPFARM~\cite{topfarm}
optimizes the objective at 8\,m/s from the same twenty unseen starts
(seeds $5$ to $24$), by sequential least-squares programming (SLSQP) on
the analytical gradients of PyWake, under the $3\,D$ separation, the
convex hull of the built positions and the two forbidden disks of the
discrete searches (exclusion zones of
TOPFARM~\cite{criadorisco2024zones}), on one process of a CPU node (\TN\
gives the configurations and the run of every seed). It moves the
machines on continuous coordinates, hence over a larger set than the
grid of the discrete searches. With the coordinates and the
objective scaled, at most $2\,000$ iterations and the continuation of
TOPFARM on the boundary (a relaxation of the constraints, as the wake
expansion continuation of~\cite{thomas2022wec} is one of the objective),
which inflates the admissible zone by
99%
~m and shrinks it to the true one
over the first 100%
\ gradient
evaluations, SLSQP converges in 7%
\
seeds and stops on a singular subproblem in the other
13%
\unskip, in every case after the
continuation has ended; its final iterate keeps the separation and the
disks in every seed and the hull in
15%
\unskip, leaving it by at most
2.64%
~m in the others. Taking
in each seed the best feasible point recorded along the run (the start
itself in 1%
\ seed), the gain is
4.30\%%
\ on average in
54%
~s, against
4.33\%%
\ for PyWake alone on the same seeds:
TOPFARM is ahead in 14%
\ of the twenty
and PyWake alone in 6%
\ (paired
difference $-$0.03%
\ points,
[$-$0.58%
, $+$0.52%
],
an interval widened by the seed at zero), and the final iterates,
feasible or not, average 4.56\%%
\unskip. On
the annual objective of
the experiment above, the plain configuration of TOPFARM (unscaled, at
most $300$ iterations) gains
+2.54\%%
\ in 391%
~s,
and the continued one +2.84\%%
\ in
315%
~s, converging in
8%
\ seeds and stopping on a singular
subproblem in the others, against
+2.94\%%
\ for PyWake alone (paired
difference $-$0.10%
\ points,
[$-$0.47%
,
$+$0.27%
]; Table~\ref{tab:proxy}). A
gradient method on continuous coordinates
thus reaches, from the same starts and on one process, gains comparable
to those of the discrete search on PyWake alone, once it is continued
and where it delivers a feasible point; what this section measures is
not the discrete search against the gradient one, but the discrete
search with and without the screening. Where the discrete variables pay
is in the constraints of Section~\ref{sec:exp-sidecon}, counts, site
costs and exclusions, which are linear on them, whereas on continuous
coordinates the counts come through a custom component or a penalty
(Section~\ref{sec:related}).
\paragraph{A second site.} To see the screening leave the site it
was built on, we repeated the comparison at 8\,m/s on Horns Rev~1: the
eighty V80 machines of the built farm ($D = 80$\,m) on their grid of
seven diameters, the twelve-sector rose of the site spread over the
degrees as for the Lillgrund instance above, the convex hull of the built
positions as the zone, dilated by half a cell diagonal so that the
machines built on its corners stay admissible, a separation of $3\,D$,
and sites on a grid at $\Delta = D$: at seven diameters of pitch the
farm is sparse, and at $D/3$ the neighborhood would exceed a million
moves. The kernel is fitted on the single wake of PyWake, on this grid
and this machine, with the recipe of Section~\ref{sec:calibration}; the start of each
seed is the built layout snapped to the grid and moved by one hundred
admissible one-exchange moves drawn with the seed, a start farther from
a local optimum than the Lillgrund ones, which are certified local
optima of the classical pairwise-Jensen guide
(Appendix~\ref{app:pt-instance}); arms, neighborhood,
counting rules and hardware are those of Table~\ref{tab:proxy}. A start
has about 136%
\ thousand one-exchange
moves, against 37\ to
79\ thousand at Lillgrund.
Table~\ref{tab:proxy} (last block) reports the outcome over
10%
\ seeds, and it differs from Lillgrund
in one respect. PyWake alone reaches a local optimum within the hour in
every seed, in 2159%
~s on average, with a
gain of +9.20\%%
\unskip. Screen-and-confirm
stops in 65%
~s at
+8.91\%%
\unskip,
0.29%
\ points below PyWake alone from
the same start (paired difference
$[$-$0.38%
,
$-$0.20%
]$), and below it in every seed.
In 6%
\ of the seeds it stops with
improving moves still in the neighborhood, which the refinement by
PyWake off the clock accepts,
27.8%
\ per seed on average in
those seeds, bringing the mean gain to
+8.99\%%
\
($-$0.21%
\ points from PyWake alone); in
the others it stops at a local optimum of PyWake lower than the one
PyWake alone reaches. The
certificate makes every run a local optimum of PyWake, in
379%
~s and at
+8.96\%%
\unskip, but these local optima
lie 0.24%
\ points
($[$-$0.34%
,
$-$0.15%
]$) below those that PyWake
alone reaches from the same starts, again in every seed. PyWake alone
takes 33.2%
\ times longer than the screen
(24.0%
\ to
50.0%
\ across the seeds) and
5.7%
\ times longer than the certified
search, in geometric mean. On this site, then, the screening keeps
its speed but not the parity of Lillgrund: the moves that the
distillate ranks first lead the search to local optima below those of
PyWake alone, by a margin small against the gain and of the same sign
in every seed. Two plausible reasons are the start, farther from a
local optimum than at Lillgrund, so that the two searches take longer
and more divergent paths, and the pitch: at seven diameters the wakes
are weaker and the moves differ by less, so that the ranking has
smaller differences to order. The mechanism is otherwise that of
Lillgrund, the distillate ranks and PyWake keeps, and the
seed-by-seed reading is in \TN.  
\section{Conclusions}
\label{sec:concl}

In this paper we have used EndoWake, introduced in our previous
paper~\cite{companion}, as a layout tool: a linear propagator distilled
once from the single-wake field of a reference model, which scores layouts
in batch on a GPU and carries zone caps, exclusions and site-dependent
costs as linear constraints on the discrete site decisions.

Its clearest benefit is computational, and it appears when EndoWake
screens moves for a common PyWake objective (Section~\ref{sec:proxy}): on
twenty starts not used to design the search, screen-and-confirm with one
GPU and thirty CPU processes reaches mean gains similar to those of the
search on PyWake alone, within the paired intervals, in a search time
about an order of magnitude shorter at one speed and on the annual
energy, and still several times shorter when the proof of local
optimality, a final scan by PyWake, is charged to the timed search.

Throughput does not, by itself, make better predictions or better
layouts. Under the RANS judge at random directions of
Section~\ref{sec:sirens-judge}, the search on PyWake is ahead on average
of the searches on EndoWake at $\Delta = D$ and at $\Delta = D/3$, but
five paired starts resolve none of the margins: the order of the means is
descriptive, and conditional on this RANS configuration. The diagnostics
of Section~\ref{sec:sirens} locate a directional quadrature error between
the directions of the rose and a spatially localized surrogate bias at
the edge of each wake, without establishing that this bias causes the
ranking. Where EndoWake is itself the objective, what the planner gains
is its form: a field that a mixed-integer model carries together with the
constraints of the project, so that an exact solver declares its gap
where a local search cannot (Section~\ref{sec:opt-tool}); at industrial
size both routes are heuristic.

These results are conditional on one site, except the comparison
of Section~\ref{sec:proxy}, repeated on a second site where the
screening kept its speedups but reached lower local optima than PyWake
alone, and those against RANS on one inflow and one speed. Refitted at other inflows, the kernel is
insensitive to the speed and sensitive to the turbulence, hence a planner
should fit it at the turbulence of the site
(Section~\ref{sec:rans-calib}); it carries none of the physics at the
scale of a farm. The RANS judge uses the published constants of its
closure, is verified on one published single wake and is compared with
the row deficits of Section~\ref{sec:rans} only; its margins are farm
power at one speed under a discrete rose, not annual energy.

Two lines of future work follow. A reduction factor that depends on the
local speed stays representable in the MILP by the device we use for the
power curve, at the price of one more piecewise-linear block per site. A
judge outside our closure, such as
TurbOPark~\cite{nygaard2020turbopark,pedersen2022turbopark} or a coarse
large-eddy simulation, would test whether the rankings of
Section~\ref{sec:sirens-judge} depend on the RANS configuration, and the
lidar wake transects of Lillgrund~\cite{sood2022lillgrund} would set the
single wake of Section~\ref{sec:rans} against a measured one. 

\section*{Funding}

The work of the first author was supported by the project
PID2024-162616OA-I00 funded by MICIU/AEI/10.13039/501100011033/FEDER, UE,
and by the Ram\'on y Cajal project RYC2023-042623-I funded by
MICIU/AEI/10.13039/501100011033 and FSE+. The funders had no role
in the study design, in the collection, analysis and interpretation of data,
in the writing of the article or in the decision to submit it for
publication.

\section*{Acknowledgments}

The RANS computations and the GPU campaigns of this paper ran on the
Upscale HPC facility of the University of Padova and on the blade cluster
of its Department of Information Engineering. We thank their staff for
the support. We are grateful to M.~Paul van der Laan and Mads Baungaard
for their comments on the RANS reference of an earlier version of this
paper, and for pointing us to their published single-wake case, against
which the closure is now verified. We also thank the two anonymous
reviewers of Applied Energy, whose reports on an earlier submission
shaped the controls and the comparisons of Sections~\ref{sec:sirens}
and~\ref{sec:proxy}.

\section*{Declaration of generative AI and AI-assisted technologies in the
manuscript preparation process}

During the preparation of this work the authors used two large language
models. Anthropic Claude served three purposes: as a programming assistant
in writing and refactoring the experimental code, as an adversarial reader
of the mathematical statements, and in editing the manuscript for language
and structure. OpenAI Codex served as an adversarial reviewer of the
manuscript. Every proposition, proof, and numerical result was verified by the
authors by executing the code, and no claim in this paper rests on a
statement the authors have not checked independently. The authors take
full responsibility for the content of the publication.

\section*{CRediT authorship contribution statement}

\textbf{Martina Fischetti}: Conceptualization, Data curation,
Funding acquisition, Investigation, Methodology, Validation, Writing -- original draft.
\textbf{Matteo Fischetti}: Conceptualization, Formal analysis,
Methodology, Software, Writing -- review \& editing.

\section*{Data availability}

The study uses no proprietary data. The power curve is the idealized cubic
curve stated in Section~\ref{sec:exp}, except where the text says that the
SWT-2.3-93 curves of PyWake are used, and on Horns Rev~1, which uses the
V80 curves of PyWake. The analytical reference is the
public Bastankhah--Port\'e-Agel formula, and the external validation uses
PyWake~2.6.20~\cite{pywake} with the
\texttt{BastankhahGaussianDeficit} model, the \texttt{LinearSum}
superposition rule and wake expansion $k=0.04$, or the value $k^{*}$ fitted
on the RANS (Table~\ref{tab:rans_calib}) in the merged-wake test of
Section~\ref{sec:exp-fidelity} and in the PyWake arm of
Section~\ref{sec:sirens}, or the fixed value
$k=0.0345$ in Section~\ref{sec:proxy}. The RANS reference
is our own OpenFOAM~\cite{openfoam} setup (Appendix~\ref{app:rans-details}).
Every instance is defined by the parameters stated with it---grid size,
cell, minimum separation, cardinality, wind rose---so that it can be rebuilt
from the text alone.

The replication package~\cite{repo} is a versioned deposit on Zenodo under
the DOI \texttt{10.5281/zenodo.22017446}, which resolves to its latest
version; it will be released publicly upon publication and is available to
the reviewers, in the meantime, through the private link given at
submission.
It holds: the instance generators and the random seeds; the OpenFOAM cases
and the disk outputs of every direction we ran, with the hub-height fields
of the single-wake, calibration and grid-convergence runs; a naive
implementation of the EndoWake sweep, written out cell by cell, which
specifies the field that the fast evaluators compute; the fast
evaluator itself, in its CUDA form and in the vectorized CPU form that
agrees with it to machine precision in double precision; the local searches, i.e.,
the vanilla local search of Section~\ref{sec:exp-heur}, the search engine
of Sections~\ref{sec:sirens} and~\ref{sec:proxy} and screen-and-confirm,
with the scripts that launch every arm; the implementation of the sector
quadrature of Section~\ref{sec:cones}; the Hexaly model of the exact
formulation; the final layouts of the searches of
Sections~\ref{sec:exp-big} to~\ref{sec:proxy}; the campaign outputs behind
the tables, which for every run of Section~\ref{sec:proxy} give the
configuration, the hardware (and, on the annual energy, the software
versions), the start, every
accepted move with its PyWake value, its time and the evaluation counts so
far, and the final layout; the rules by which times and evaluations are
counted; the raw timing logs of Section~\ref{sec:exp-batch}; the scripts
that write each table directly from the campaign outputs; and a technical
note, generated from the same outputs, with the fluid-dynamic properties of
EndoWake, the validation instrumentation, the record of the exact
optimization and of the RANS reference, and the GPU implementation notes.
Our previous paper~\cite{companion} has its own package. With the
package, a reader can rerun every search of Section~\ref{sec:proxy} on a
node with one GPU, and recompute, with PyWake alone, the value of every
start, accepted move and final layout from the logs, and every verdict of
the RANS judge from the OpenFOAM cases.

\section*{Declaration of competing interest}

The authors declare that they have no known competing financial
interests or personal relationships that could have appeared to influence
the work reported in this paper. 

\appendix

\counterwithin{equation}{section}
\counterwithin{table}{section}
\counterwithin{figure}{section}

\section{The RANS reference in detail}
\label{app:rans-details}

In this appendix we complement Table~\ref{tab:rans_setup} with the
verification of the closure, the average over the wind direction, the
monitor, the controls of grid and noise, and the random judge of
Section~\ref{sec:rans-test}. The reference remains a RANS with its own
closure and grid and its own errors, not a measurement of how a layout
performs at the site; \TN\ holds the record of the runs and the tables
summarized here.

\paragraph{Verification of the closure.} OpenFOAM does not ship the
$k$-$\varepsilon$-$f_P$ model of van der Laan et
al.~\cite{vanderlaan2015kefp,vanderlaan2015double,vanderlaan2015farms},
so we implemented it as a RANS model derived from the standard
$k$-$\varepsilon$ class, with eddy viscosity $\nu_t = C_\mu f_P
k^2/\varepsilon$, where
\[
f_P \;=\; \frac{2 f_0}{1 + \sqrt{1 + 4 f_0 (f_0 - 1)\,(\sigma/\tilde\sigma)^2}},
\qquad
f_0 = \frac{C_R}{C_R - 1},
\qquad
\sigma = \frac{k}{\varepsilon}\,|\nabla U|,
\qquad
\tilde\sigma = \frac{1}{\sqrt{C_\mu}} .
\]
The production term of the $\varepsilon$ equation is written with the
same eddy viscosity, $f_P$ included. Upstream of the farm $f_P$ equals one
to three decimal places; one to two diameters behind the first machine it
stays near one outside the wake and on its axis, and drops well below one
in the shear layer---the shape of the published contours. As a
quantitative check we ran the single-turbine case of~\cite{baungaard2022}
in the version distributed with~\cite{actuatordiskfoam} ($D = 80$~m,
thrust coefficient $0.75$, 8\,m/s, turbulence intensity $5.7\%$), which
also distributes the hub-height field of the authors' solver, with the
same closure and constants. Along
the wake axis, over the stations between $3$ and $10\,D$, our speed
differs from theirs by 0.033\ of the free
stream in root mean square (0.72\
against 0.69\ at $5\,D$), and by
0.048\ when the speed is averaged
over the rotor width. The constants are those of the
authors of the closure. The previous submission of this work used
$C_R = 60$: in OpenFOAM v2412 the production term of the
$\varepsilon$ equation of the base class is written for the eddy
viscosity of the standard model, and the implementation of that
version, which redefined the eddy viscosity only, left $f_P$ out of
that term; the wake then recovered too slowly at $C_R = 4.5$, and that
version compensated with $C_R = 60$, chosen to match the row deficits
of the site. With the term corrected, $C_R = 4.5$ reproduces the
published single wake to the accuracy stated above, and the calibration,
the judgments and the controls of this paper were recomputed with the
corrected closure; the searches of Section~\ref{sec:proxy}, whose
objective is PyWake, keep a fixed $k$, as that section states. \TN\
records the corrected closure and its verification.

\paragraph{The inflow along an empty domain.} The rough-wall condition on
$\varepsilon$ keeps the default constants of OpenFOAM ($C_\mu = 0.09$,
$\kappa = 0.41$), while the inlet profiles use those of the closure
(Table~\ref{tab:rans_setup}), so the profile is not exactly in
equilibrium. We measured its drift on the case of the judge at
$222^\circ$ with no thrust on the disks (the setting $C_T = 0.001$, at
which the momentum method of the solver gives zero axial induction), the
same grid and the same iterations. Averaged across the width of the
domain at hub height, in bands of $25$~m along the wind from the first
machine to the last, the speed exceeds its inlet value by at most
1.2\%\ (0.7\%\
on average), and the turbulence intensity falls below the inlet value
of Table~\ref{tab:rans_setup} by at most
5.4\%\ of that value
(3.3\%\ on average). The drift is a
gradient along the wind, the same in every judged run; the gains,
measured over the start at the same directions, largely cancel it.

\paragraph{The average over the wind direction.} We ran the built layout
at $222^\circ + \delta$ and at $120^\circ + \delta$, for $\delta = 0, \pm
1, \pm 2, \pm 4, \pm 6, \pm 8$ degrees, on the grid and with the monitor
of the judge, and report in Section~\ref{sec:rans} the power of the second
turbine of the aligned row relative to the first at the single direction
and averaged with a Gaussian weight of standard deviation $3.3^\circ$,
the uncertainty of the wind direction estimated for this farm
in~\cite{vanderlaan2015farms}, truncated at $\pm 8^\circ$ and
renormalized; the uniform average over a bin of $\pm 2.5^\circ$ lies
between the two. The measured values are bin averages that also carry
this uncertainty~\cite{gaumond2014}.

\paragraph{The monitor.} The reference velocity of each actuator disk is
the mean of nine points over the rotor, one diameter upstream. Against a
single point on the axis, the nine points see more wind at the rotor
edges and raise the second-turbine ratio by
0.06\ with standard
$k$-$\varepsilon$. The nine points are probes averaged by the solver: in
OpenFOAM v2412 a cell-set monitor would read the velocity inside the
rotor itself, and under-thrust.

\paragraph{Grid convergence.} On the built layout at $222^\circ$, with
$k$-$\varepsilon$ at the stock OpenFOAM constants ($C_\mu = 0.09$) and the
axis monitor, halving every cell changes the farm power by
+2.0\%\ and the second-turbine ratio
(0.379 $\to$ 0.395\unskip), and leaves the rows
from the fourth on unchanged: a shift of level. With the judging closure
and monitor, the second-turbine ratio of the same row moves from
0.23\ to
0.24\ on the halved grid. Re-judged
on the halved grid at the six sector centers, for seeds~$0$ and~$1$ of
the bench of Section~\ref{sec:rans-test}, the margin between arms~A
and~B changes by +0.10\ points, whereas
the level moves with the grid.

\paragraph{Noise of the judge.} On the common start and the final layouts
of arms~A and~B of seed~$0$, under the judging closure at $120^\circ$,
the aligned direction with the tightest spacing, two identical runs at
the tolerance of the calibration runs differ by at most
0.03\%\ in farm power, and the tenfold tighter
tolerance of the judged runs moves the level by at most
0.04\%\unskip. Moving one turbine by half a
diameter changes the farm power by an amount that the two tolerances
reproduce to within 0.05\ percentage
points, and that the halved base cell reproduces to within
0.03\ points: on these layouts and at
this direction, the coarse grid resolves a half-diameter move. The judged
run of arm~B ends at iteration 1\,500\unskip,
with initial residuals of $3.8 \times 10^{-8}$\ on the
streamwise velocity and $2.2 \times 10^{-6}$\ on the
pressure, and the farm power drifts by $<$0.001\%\
over the last $200$ iterations. Over all the
552\ judged runs counted in
Table~\ref{tab:rans_setup}, starts included (the
240%
\ runs of the controlled comparison
of Section~\ref{sec:sirens-judge} are audited separately in \TN), the
final residuals are at
most $1.9\cdot10^{-7}$\ on the streamwise
velocity and $3.9\cdot10^{-6}$\ on the
pressure, within the control in every case; those of $k$ and
$\varepsilon$ end between
$1.8\cdot10^{-6}$\ (median) and
$2.3\cdot10^{-6}$\unskip, so that the
solver reports convergence on all the fields in
24\% of the runs
(25\% end with the
recorded residuals of $k$ and $\varepsilon$ within the thresholds) and
the others stop at the cap of $1\,500$ iterations. Over the last five iterations the farm power of every run
moves by at most $7.2\cdot10^{-5}$\%\ of
its value, and in the 8\
repeated runs at the tolerance of the judged runs, by at most
$1.5\cdot10^{-3}$\%\ over the
last $200$ iterations. The 6%
\ capped
runs with the largest residuals, rerun with the cap doubled to
3\,000%
\ iterations, change their
farm power by at most $8.7\cdot10^{-6}$\%%
\ beyond
iteration $1\,500$ (\TN). These controls bound the numerical noise of the judgments
they cover, not whether the reference ranks layouts as the site would.

\paragraph{The random judge.} The directions of
Section~\ref{sec:rans-test} are drawn once from a fixed seed, one
uniformly in each thirty-degree interval between consecutive sector
centers (\TN\ lists them). Since the intervals run between the sector
centers, whereas the sectors are centered on them, the draw is not a
stratified quadrature over the sectors: a sector may hold zero, one or
two directions. Of the 60\ drawn
directions, 9\ are among the thirty
directions, twelve degrees apart, of the rose of arm~A, and none is a
Gauss--Legendre node of arm~B. 
\section{Conventions, environment and numerical hygiene}
\label{app:exp-conventions}

\paragraph{Environment.} Table~\ref{tab:env} reports where each campaign
ran. The laptop GPU (GTX 1050 Ti\unskip) is the slow
end of the hardware a reader is likely to have, measured on mains power,
where the card holds its boost clock under load (median
1670~MHz). We implemented the model in
Python~3.14 through \texttt{gurobipy} and used Gurobi~13.0 as MILP solver
and Hexaly~15.0 as heuristic solver. The power curve enters the MILP
through the SOS2 block~\eqref{eq:mip-power}, with the breakpoints of
Section~\ref{sec:mip}, throughout. The budget per solve, the thread count
and the number of seeds of each campaign are stated in the caption of its
table.

\paragraph{The proxy search.} Each seed of Section~\ref{sec:proxy} ran on
one node with an AMD EPYC~9335 CPU (32 cores) and an NVIDIA RTX PRO~6000
Blackwell Server Edition. The runs on the annual energy logged the
software versions: PyWake~2.6.20, CuPy~14.2.0, NumPy~2.5.2,
Python~3.13.5,
SciPy~1.17.1, the CUDA runtime
13.2\ and a driver supporting CUDA
13.4\unskip; the other runs of the section
logged the node, the processes and the GPU, not the versions. Its times are \emph{warmed search wall times}, from the
evaluation of the start to the end of the search. They include the
enumeration of every neighborhood and the calibration of the list, and
exclude the warm-up of both models, the compilation of the GPU kernels,
the start of the process pool and the distillation, which a first search
pays on top.

\begin{table}[htbp]
\centering
\caption{Where each campaign ran. The cluster GPUs are used in double
precision, except in the proxy search and in the single-precision batch of
Section~\ref{sec:exp-batch}; the laptop GPU is the floor of the throughput figures.}
\label{tab:env}
\begin{tabular}{p{0.40\textwidth}p{0.54\textwidth}}
\toprule
campaign & machine \\
\midrule
MILP and local search (Section~\ref{sec:opt-tool}) &
laptop, Intel Core i9-8950HK at $2.9$\,GHz, $32$\,GB of RAM, Windows~11 \\
batched evaluator, one direction (Section~\ref{sec:exp-batch}) &
NVIDIA RTX PRO 6000 Blackwell Server Edition\ of the university cluster; the laptop GPU
(GTX 1050 Ti\unskip) alongside as the floor \\
batched evaluator, twelve-sector rose (Section~\ref{sec:exp-batch}) &
laptop GPU only \\
RANS reference and judge (Sections~\ref{sec:rans} and~\ref{sec:rans-test}) &
CPU nodes of the university cluster, OpenFOAM \\
RANS-judged campaigns on Lillgrund
(Section~\ref{sec:sirens}) &
PyWake arm on one CPU node, thirty worker processes; endogenous arms on
data-center GPUs (NVIDIA RTX PRO 6000) \\
proxy search (Section~\ref{sec:proxy}) &
one GPU node per seed (AMD EPYC~9335, 32 cores), all arms in sequence on
it: PyWake on 30\ CPU processes in every arm
(eight in the repeat on eight processes), EndoWake on one RTX PRO 6000 in
single precision \\
\bottomrule
\end{tabular}
\end{table}

\paragraph{Numerical hygiene.} We set Gurobi's \texttt{NumericFocus} to
its highest value throughout, as a precaution: in the big-$M$ wake
implication~\eqref{eq:mip-wake} a binary within the integrality
tolerance leaves a slack that the constant multiplies and the
substitution chains carry downstream, so that a certificate at default
tolerances is only as tight as that slack. On the largest instance that
both settings close, the setting multiplies solution time by
1.5$\times$. EndoWake runs in single precision in
the throughput measurements of Section~\ref{sec:exp-batch} and in the
proxy search, where it only scores and ranks; the searches of
Section~\ref{sec:sirens} ran in double precision, which the accuracy of
the model does not require but which keeps every run identical across
cards under the same version of the array library (\TN). 

\section{Instances and protocols}
\label{app:protocols}

In this appendix we give the instances and protocols of the two search
campaigns on Lillgrund, one per subsection; \TN\ tabulates, for every
search campaign of Sections~\ref{sec:opt-tool} to~\ref{sec:proxy}, the
power curve, the objective, the judge and the variant of the EndoWake
field.

\subsection{The bench of Section~\ref{sec:rans-test}}
\label{app:pt-instance}

\paragraph{The instance.} The grid is aligned with the UTM axes, rows
running south and columns east, with $D = 93$\,m;
its north-west corner lies $2\,D$ west and north of the built farm, at
easting 358\,619\,m and northing
6\,155\,495\,m. At $\Delta = D$ the grid is
$34 \times 33$ and the sites are the
529\ cells inside the convex hull of the
built positions. Around the two lattice sites left empty inside the hull
we exclude disks of radius 1.65\,$D$ (half
the smaller pitch), whose cells are unavailable but are not obstacles to
the wind. The field of arm~B$'$ has
$100 \times 98$ cells at $\Delta = D/3$, one in
three in each direction allocatable, i.e., the same sites.

\paragraph{The start.} The guide that builds the start of each seed scores
a layout by pairwise Jensen deficits with $k = 0.05$ in free stream,
summed linearly, under the idealized curve. Its greedy construction draws a
shortlist half from the windiest free cells of the rose-weighted field and
half at random, and takes the best candidate on the full objective; a
one-exchange search then reaches a certified local optimum. Under a cap of
$48$ machines, the construction places what the separation admits, namely
45/47/47/45/46\ for seeds $0$--$4$.

\paragraph{The arms.}\label{app:pt-protocol}
Arm~B and the guide read the twelve sectors of the Lillgrund dataset of
PyWake~\cite{pywake} (the sector at $60^\circ$ lies below cut-in and
contributes nothing), each at its Weibull mean speed and frequency;
integrated by five\ Gauss--Legendre nodes per
sector, they give 60\ scenarios. Arms~B and
B$'$ carry the kernel of Table~\ref{tab:rans_calib} at the cell of their
field. Arm~A runs PyWake~2.6.20 on the SWT-2.3-93 of the same dataset, at
30\ directions
12$^\circ$ apart and speeds from $3$ to
$25$\,m/s, with \texttt{Bastankhah\_PorteAgel\_2014} at its defaults
(\texttt{LinearSum} superposition, speed sampled at the rotor center,
deficit scaled by the local speed at the upstream rotor, no turbulence or
deflection model) except for $k = k^{*}$. Note that $k^{*}$ was fitted with
the thrust coefficient $8/9$ of the
idealized machine, whereas PyWake applies it with the thrust curve of the
SWT-2.3-93; refitted with the thrust of the SWT-2.3-93 at $8$\,m/s, it
changes the objective of arm~A by at most
0.099\ points on the start and final
layouts of the five seeds.

\paragraph{The search.} Arm~A visits the turbines in random order and
moves each to the best of 256\ sampled
admissible cells if this improves the objective; a run stops after
eight\ consecutive sweeps without improvement,
on 30\ worker processes. Arm~B scans in
blocks of 8\,192\ layouts. At $\Delta = D$, arm~A
restarted 53\ times per seed on average
(48--58) and arm~B
1\,115\
(947--1\,483).

\subsection{The proxy search of Section~\ref{sec:proxy}}
\label{app:sm-proxy}

\emph{The fit.} We fit the kernel by the recipe of
Section~\ref{sec:calibration}, with the window from three to twelve
diameters and the footprint of~\eqref{eq:mip-wake}, on the hub-height field
of one isolated SWT-2.3-93 computed by PyWake~2.6.20
(\texttt{Bastankhah\_PorteAgel\_2014} at $k = 0.0345$, 8\,m/s) on the grid at
$\Delta = D/3$. It returns $\alpha = 0.083$,
$\beta = 0$, $\gamma = 0.0345$ and a source
amplitude of 0.644\unskip, with a residual of
6.8\%\unskip.

\emph{The neighborhood and the arms.} A move of Section~\ref{sec:proxy} is
improving if it raises the value of the incumbent by more than $10^{-9}$ in
the units of the evaluator, as in Section~\ref{sec:rans-test}. Every arm
draws its random order from the generator of the seed, so that the samples
of screen-and-confirm are nested. The block size of PyWake alone
starts at the number of processes, doubles after every block without an
improving move (up to $2\,048$), and returns to its minimum after an
accepted move. EndoWake alone runs the same scan on the distillate, in
single precision. In screen-and-confirm, PyWake evaluates the list in
blocks of as many moves as it has processes: within a block the first
improving move in rank order is accepted, and the rest of the block
returns to the head of the list to be evaluated again on the new
incumbent. The rate measurement draws its neighbors of the start from a
generator of its own and charges its evaluations to the search. The fixed
list of Section~\ref{sec:proxy} confirms $\lceil 3\,000\, q / |\mathcal{N}|
\rceil$ moves of a sample of $q$ out of a neighborhood of $|\mathcal{N}|$.
The run of PyWake alone that checks a final layout as a local optimum
starts from that layout, off the clock, and runs to a local optimum of
PyWake. Each seed is one job on one GPU node, with the arms in sequence
(Table~\ref{tab:env}); for seed~$5$, we re-evaluated the trajectory of
every arm with PyWake through an independent code path, and found
agreement up to rounding.

\emph{The annual energy.} The protocol was fixed in writing before any
search on the annual energy: the seeds, with twenty reserve seeds to be run
only if the half-width of the 95\% interval of the paired difference of the
full screen exceeded 0.10\ points (it was
0.08); four arms in a fixed order on one node in
the original campaign (PyWake alone, the full screen, the full screen
with the certificate, the reduced screen), and the reduced screen
with the certificate in a later campaign on the same starts, node type
and frozen protocol (\TN); a limit of 7\,200~s
per search, which no run reached; and the checks below. Each direction carries a
thirtieth of the frequency of its sector; power and thrust are interpolated
linearly in the tables of PyWake, with cut-in at 3\,m/s, cut-out at 25\,m/s
and no thrust when idle, in both models. Each node of the reduced screen
carries, direction by direction, the Weibull mass of the bins from 3 to
7\,m/s, 7 to 9, 9 to 11 and 11 to 13 respectively, and uses the fit at its
own speed. EndoWake ranks in blocks of 8\,192\
layouts, the size of highest rate on the full objective in a benchmark run
before the searches. Every final layout of the two screens was checked off
the clock against its whole neighborhood
(1\,131~s on average), without moving it, and
the start and the final layout of every run were evaluated again with the
step halved in direction, in speed and in both, and with the rose rotated by
half a degree; these re-evaluations keep the layouts fixed, and the paired
difference of a single seed moves by at most
0.026\ points. Charging each run with its own
setup (loading, the process pool, the construction and the warm-up of both
models) moves the speedup of the reduced screen to
14.5; the fit of EndoWake
(7~s) is shared by all the runs. The seeds
ran on two groups of nodes with the same processor, GPU and libraries,
which differ in the minor version of the operating-system kernel and in
the build of Python; the four arms of a seed share one node, and the
speedup of the reduced screen is 12.6\ on
one group and 17.7\ on the other.


\begin{thebibliography}{99}
\bibitem{fischetti2019overview}
M.~Fischetti, D.~Pisinger.
\newblock Mathematical optimization and algorithms for offshore wind farm
design: an overview.
\newblock {\em Business \& Information Systems Engineering},
61(4):469--485, 2019.
\newblock \texttt{doi:10.1007/s12599-018-0538-0}.

\bibitem{fischetti2020vattenfall}
M.~Fischetti, J.~R. Kristoffersen, T.~Hjort, M.~Monaci, D.~Pisinger.
\newblock Vattenfall optimizes offshore wind farm design.
\newblock {\em INFORMS Journal on Applied Analytics}, 50(1):80--94, 2020.
\newblock \texttt{doi:10.1287/inte.2019.1019}.

\bibitem{jensen1983}
N.~O. Jensen.
\newblock A note on wind generator interaction.
\newblock Technical Report Ris\o-M-2411, Ris\o\ National Laboratory, Roskilde,
1983.

\bibitem{bastankhah2014}
M.~Bastankhah, F.~Port\'e-Agel.
\newblock A new analytical model for wind-turbine wakes.
\newblock {\em Renewable Energy}, 70:116--123, 2014.
\newblock \texttt{doi:10.1016/j.renene.2014.01.002}.

\bibitem{pena2018anholt}
A.~Pe\~na, K.~S. Hansen, S.~Ott, M.~P. van~der Laan.
\newblock On wake modeling, wind-farm gradients, and AEP predictions at the
Anholt wind farm.
\newblock {\em Wind Energy Science}, 3(1):191--202, 2018.
\newblock \texttt{doi:10.5194/wes-3-191-2018}.

\bibitem{sebastiani2021}
A.~Sebastiani, F.~Castellani, G.~Crasto, A.~Segalini.
\newblock Data analysis and simulation of the Lillgrund wind farm.
\newblock {\em Wind Energy}, 24(6):634--648, 2021.
\newblock \texttt{doi:10.1002/we.2594}.

\bibitem{king2017}
R.~N. King, K.~Dykes, P.~Graf, P.~E. Hamlington.
\newblock Optimization of wind plant layouts using an adjoint approach.
\newblock {\em Wind Energy Science}, 2(1):115--131, 2017.
\newblock \texttt{doi:10.5194/wes-2-115-2017}.

\bibitem{antonini2020adjoint}
E.~G.~A. Antonini, D.~A. Romero, C.~H. Amon.
\newblock Optimal design of wind farms in complex terrains using computational
fluid dynamics and adjoint methods.
\newblock {\em Applied Energy}, 261:114426, 2020.
\newblock \texttt{doi:10.1016/j.apenergy.2019.114426}.

\bibitem{turner2014}
S.~D.~O. Turner, D.~A. Romero, P.~Y. Zhang, C.~H. Amon, T.~C.~Y. Chan.
\newblock A new mathematical programming approach to optimize wind farm
layouts.
\newblock {\em Renewable Energy}, 63:674--680, 2014.
\newblock \texttt{doi:10.1016/j.renene.2013.10.023}.

\bibitem{zhang2014}
P.~Y. Zhang, D.~A. Romero, J.~C. Beck, C.~H. Amon.
\newblock Solving wind farm layout optimization with mixed integer programs
and constraint programs.
\newblock {\em EURO Journal on Computational Optimization}, 2(3):195--219,
2014.
\newblock \texttt{doi:10.1007/s13675-014-0024-5}.

\bibitem{fischetti2016proximity}
M.~Fischetti, M.~Monaci.
\newblock Proximity search heuristics for wind farm optimal layout.
\newblock {\em Journal of Heuristics}, 22(4):459--474, 2016.
\newblock \texttt{doi:10.1007/s10732-015-9283-4}.

\bibitem{fischetti2021nextgen}
M.~Fischetti.
\newblock On the optimized design of next-generation wind farms.
\newblock {\em European Journal of Operational Research},
291(3):862--870, 2021.
\newblock \texttt{doi:10.1016/j.ejor.2020.10.048}.

\bibitem{fischetti2023integrated}
M.~Fischetti, M.~Fischetti.
\newblock Integrated layout and cable routing in wind farm optimal design.
\newblock {\em Management Science}, 69(4):2147--2164, 2023.
\newblock \texttt{doi:10.1287/mnsc.2022.4470}.

\bibitem{perezrua2023}
J.-A. P\'erez-R\'ua, M.~Stolpe, N.~A. Cutululis.
\newblock A neighborhood search integer programming approach for wind farm
layout optimization.
\newblock {\em Wind Energy Science}, 8(9):1453--1473, 2023.
\newblock \texttt{doi:10.5194/wes-8-1453-2023}.

\bibitem{cazzaro2022vns}
D.~Cazzaro, D.~Pisinger.
\newblock Variable neighborhood search for large offshore wind farm layout
optimization.
\newblock {\em Computers \& Operations Research}, 138:105588, 2022.
\newblock \texttt{doi:10.1016/j.cor.2021.105588}.

\bibitem{cazzaro2022multiscale}
D.~Cazzaro, A.~Trivella, F.~Corman, D.~Pisinger.
\newblock Multi-scale optimization of the design of offshore wind farms.
\newblock {\em Applied Energy}, 314:118830, 2022.
\newblock \texttt{doi:10.1016/j.apenergy.2022.118830}.

\bibitem{topfarm}
DTU Wind Energy.
\newblock TOPFARM: a wind plant optimization framework, version 2 [software].
\newblock DTU Wind, Technical University of Denmark, 2026.
\newblock \texttt{https://topfarm.pages.windenergy.dtu.dk/TopFarm2/} (accessed
September 2026).

\bibitem{pywake}
M.~M. Pedersen, A.~Meyer Forsting, P.~van~der Laan, R.~Riva,
L.~A. Alcayaga Rom\`an, J.~Criado Risco, et al.
\newblock PyWake: an open-source wind farm simulation tool, version 2.6.20
[software].
\newblock DTU Wind, Technical University of Denmark, 2026 (release of 30 March
2026).
\newblock \texttt{https://gitlab.windenergy.dtu.dk/TOPFARM/PyWake} (accessed
September 2026).

\bibitem{floris}
NREL.
\newblock {\em FLORIS: FLOw Redirection and Induction in Steady State},
version 4 [software].
\newblock National Renewable Energy Laboratory, 2024.
\newblock \texttt{https://github.com/NREL/floris} (accessed September 2026).

\bibitem{locascio2022}
M.~J.~LoCascio, C.~J.~Bay, M.~Bastankhah, G.~E.~Barter, P.~A.~Fleming,
L.~A.~Mart\'inez-Tossas.
\newblock {FLOW} Estimation and Rose Superposition ({FLOWERS}): an integral
approach to engineering wake models.
\newblock {\em Wind Energy Science}, 7(3):1137--1151, 2022.
\newblock \texttt{doi:10.5194/wes-7-1137-2022}.

\bibitem{locascio2024}
M.~J.~LoCascio, C.~J.~Bay, L.~A.~Mart\'inez-Tossas, M.~Bastankhah, C.~Gorl\'e.
\newblock {FLOWERS AEP}: an analytical model for wind farm layout
optimization.
\newblock {\em Wind Energy}, 27(12):1563--1580, 2024.
\newblock \texttt{doi:10.1002/we.2954}.

\bibitem{locascio2025}
M.~J.~LoCascio, C.~J.~Bay, L.~A.~Mart\'inez-Tossas, J.~J.~Thomas,
C.~Gorl\'e.
\newblock Efficient wind farm layout optimization with the {FLOWERS AEP}
model and analytic gradients.
\newblock {\em Journal of Renewable and Sustainable Energy}, 17(1):013302,
2025.
\newblock \texttt{doi:10.1063/5.0237778}.

\bibitem{criadorisco2024zones}
J.~Criado Risco, R.~Valotta Rodrigues, M.~Friis-M{\o}ller, J.~Quick,
M.~M{\o}lgaard Pedersen, P.-E. R\'ethor\'e.
\newblock Gradient-based wind farm layout optimization with inclusion and
exclusion zones.
\newblock {\em Wind Energy Science}, 9(3):585--600, 2024.
\newblock \texttt{doi:10.5194/wes-9-585-2024}.

\bibitem{companion}
M.~Fischetti, M.~Fischetti.
\newblock EndoWake: Modeling wind by linear programming for wind farm layout
optimization.
\newblock Preprint, \texttt{arXiv:2609.40157}, 2026; submitted
to an international journal.
\newblock Replication package on Zenodo,
\texttt{doi:10.5281/zenodo.22017424}.

\bibitem{zong2020}
H.~Zong, F.~Port\'e-Agel.
\newblock A momentum-conserving wake superposition method for wind farm power
prediction.
\newblock {\em Journal of Fluid Mechanics}, 889:A8, 2020.
\newblock \texttt{doi:10.1017/jfm.2020.77}.

\bibitem{lanzilao2022}
L.~Lanzilao, J.~Meyers.
\newblock A new wake-merging method for wind-farm power prediction in the
presence of heterogeneous background velocity fields.
\newblock {\em Wind Energy}, 25(2):237--259, 2022.
\newblock \texttt{doi:10.1002/we.2669}.

\bibitem{bastankhah2021}
M.~Bastankhah, B.~L.~Welch, L.~A.~Mart\'inez-Tossas, J.~King, P.~Fleming.
\newblock Analytical solution for the cumulative wake of wind turbines in wind
farms.
\newblock {\em Journal of Fluid Mechanics}, 911:A53, 2021.
\newblock \texttt{doi:10.1017/jfm.2020.1037}.

\bibitem{zhang2026}
Z.~Zhang, J.~He, Y.~Liu, Q.~Li, A.~Kareem.
\newblock A new physically based model for cumulative wakes and application to
wind farms.
\newblock {\em Energy Conversion and Management}, 361:121593, 2026.
\newblock \texttt{doi:10.1016/j.enconman.2026.121593}.

\bibitem{smith2006}
J.~E. Smith, R.~L. Winkler.
\newblock The optimizer's curse: skepticism and postdecision surprise in
  decision analysis.
\newblock {\em Management Science}, 52(3):311--322, 2006.
\newblock \texttt{doi:10.1287/mnsc.1050.0451}.

\bibitem{thomas2023les}
J.~J. Thomas, C.~J. Bay, A.~P.~J. Stanley, A.~Ning.
\newblock Gradient-based wind farm layout optimization results compared
with large-eddy simulations.
\newblock Preprint, {\em Wind Energy Science Discussions}, 2022; the
revision was not accepted.
\newblock \texttt{doi:10.5194/wes-2022-4}.

\bibitem{schoeler2026validation}
J.~P. Sch{\o}ler, E.~Simutis, M.~P. van~der Laan, J.~Quick, P.-E. R{\'e}thor{\'e}.
\newblock Validation of {RANS}-calibrated engineering models and {ANN}-based
surrogate for wind farm flow simulation and layout optimization.
\newblock {\em Wind Energy Science}, 11:2845--2868, 2026.
\newblock \texttt{doi:10.5194/wes-11-2845-2026}.

\bibitem{fischetti2017thesis}
M.~Fischetti.
\newblock Mathematical programming models and algorithms for offshore wind
park design.
\newblock PhD thesis, DTU Management Engineering, Technical University of
Denmark, Kgs.~Lyngby, 2017.

\bibitem{feng2015}
J.~Feng, W.~Z. Shen.
\newblock Modelling wind for wind farm layout optimization using joint
distribution of wind speed and wind direction.
\newblock {\em Energies}, 8(4):3075--3092, 2015.
\newblock \texttt{doi:10.3390/en8043075}.

\bibitem{quaeghebeur2021}
E.~Quaeghebeur, R.~Bos, M.~B.~Zaaijer.
\newblock Wind farm layout optimization using pseudo-gradients.
\newblock {\em Wind Energy Science}, 6(3):815--839, 2021.
\newblock \texttt{doi:10.5194/wes-6-815-2021}.

\bibitem{jasa2022multifidelity}
J.~Jasa, P.~Bortolotti, D.~Zalkind, G.~Barter.
\newblock Effectively using multifidelity optimization for wind turbine design.
\newblock {\em Wind Energy Science}, 7(3):991--1006, 2022.
\newblock \texttt{doi:10.5194/wes-7-991-2022}.

\bibitem{repo}
M.~Fischetti, M.~Fischetti.
\newblock Replication package and technical note for ``GPU-accelerated wind
farm layout search with a distilled endogenous wake model (EndoWake)''
[dataset].
\newblock Zenodo, 2026. \texttt{doi:10.5281/zenodo.22017446} (all versions);
to be released publicly at publication, and available to the reviewers
through the private link given at submission.

\bibitem{burton2011}
T.~Burton, N.~Jenkins, D.~Sharpe, E.~Bossanyi.
\newblock {\em Wind Energy Handbook}, 2nd edition.
\newblock Wiley, 2011.
\newblock \texttt{doi:10.1002/9781119992714}.

\bibitem{porteagel2020}
F.~Port\'e-Agel, M.~Bastankhah, S.~Shamsoddin.
\newblock Wind-turbine and wind-farm flows: a review.
\newblock {\em Boundary-Layer Meteorology}, 174(1):1--59, 2020.
\newblock \texttt{doi:10.1007/s10546-019-00473-0}.

\bibitem{lissaman1979}
P.~B.~S. Lissaman.
\newblock Energy effectiveness of arbitrary arrays of wind turbines.
\newblock {\em Journal of Energy}, 3(6):323--328, 1979.
\newblock \texttt{doi:10.2514/3.62441}.

\bibitem{katic1986}
I.~Kati\'c, J.~H\o jstrup, N.~O. Jensen.
\newblock A simple model for cluster efficiency.
\newblock In W.~Palz, E.~Sesto, editors, {\em EWEC'86: Proceedings of the
European Wind Energy Association Conference and Exhibition, Rome, 1986},
volume~1, pages 407--410. A.~Raguzzi, Rome, 1987.

\bibitem{niayifar2016}
A.~Niayifar, F.~Port\'e-Agel.
\newblock Analytical modeling of wind farms: a new approach for power
prediction.
\newblock {\em Energies}, 9(9):741, 2016.
\newblock \texttt{doi:10.3390/en9090741}.

\bibitem{locascio2023}
M.~J.~LoCascio, C.~Gorl\'e, M.~F.~Howland.
\newblock Data-driven wake model parameter estimation to analyze effects of
wake superposition.
\newblock {\em Journal of Renewable and Sustainable Energy}, 15(6):063304,
2023.
\newblock \texttt{doi:10.1063/5.0163896}.

\bibitem{anagnostopoulos2023wakenet}
S.~J. Anagnostopoulos, J.~Bauer, M.~C.~A. Clare, M.~D. Piggott.
\newblock Accelerated wind farm yaw and layout optimisation with multi-fidelity
deep transfer learning wake models.
\newblock {\em Renewable Energy}, 218:119293, 2023.
\newblock \texttt{doi:10.1016/j.renene.2023.119293}.

\bibitem{padron2019pce}
A.~S. Padr\'on, J.~Thomas, A.~P.~J. Stanley, J.~J. Alonso, A.~Ning.
\newblock Polynomial chaos to efficiently compute the annual energy production
in wind farm layout optimization.
\newblock {\em Wind Energy Science}, 4(2):211--231, 2019.
\newblock \texttt{doi:10.5194/wes-4-211-2019}.

\bibitem{quick2023sgd}
J.~Quick, P.-E. R\'ethor\'e, M.~M{\o}lgaard Pedersen, R.~Valotta Rodrigues,
M.~Friis-M{\o}ller.
\newblock Stochastic gradient descent for wind farm optimization.
\newblock {\em Wind Energy Science}, 8(8):1235--1250, 2023.
\newblock \texttt{doi:10.5194/wes-8-1235-2023}.

\bibitem{valottarodrigues2024speeding}
R.~Valotta Rodrigues, M.~M{\o}lgaard Pedersen, J.~P. Sch{\o}ler, J.~Quick,
P.-E. R\'ethor\'e.
\newblock Speeding up large-wind-farm layout optimization using gradients,
parallelization, and a heuristic algorithm for the initial layout.
\newblock {\em Wind Energy Science}, 9(2):321--341, 2024.
\newblock \texttt{doi:10.5194/wes-9-321-2024}.

\bibitem{donovan2005}
S.~Donovan.
\newblock Wind farm optimization.
\newblock In {\em Proceedings of the 40th Annual Conference of the Operational
Research Society of New Zealand}, pages 196--205, Wellington, 2005.

\bibitem{fagerfjall2010}
P.~Fagerfj\"all.
\newblock Optimizing wind farm layout---more bang for the buck using
mixed integer linear programming.
\newblock Master's thesis, Chalmers University of Technology and University of
Gothenburg, 2010.

\bibitem{kuo2016}
J.~Y.~J. Kuo, D.~A. Romero, J.~C. Beck, C.~H. Amon.
\newblock Wind farm layout optimization on complex terrains---integrating a
CFD wake model with mixed-integer programming.
\newblock {\em Applied Energy}, 178:404--414, 2016.
\newblock \texttt{doi:10.1016/j.apenergy.2016.06.085}.

\bibitem{thomas2023}
J.~J. Thomas, N.~F. Baker, P.~Malisani, E.~Quaeghebeur, S.~Sanchez
Perez-Moreno, J.~Jasa, et al.
\newblock A comparison of eight optimization methods applied to a wind farm
layout optimization problem.
\newblock {\em Wind Energy Science}, 8(5):865--891, 2023.
\newblock \texttt{doi:10.5194/wes-8-865-2023}.

\bibitem{baker2019}
N.~F. Baker, A.~P.~J. Stanley, J.~J. Thomas, A.~Ning, K.~Dykes.
\newblock Best practices for wake model and optimization algorithm selection
in wind farm layout optimization.
\newblock In {\em AIAA Scitech 2019 Forum}, San Diego, CA, 2019. Paper AIAA
2019-0540.
\newblock \texttt{doi:10.2514/6.2019-0540}.

\bibitem{guirguis2016}
D.~Guirguis, D.~A. Romero, C.~H. Amon.
\newblock Toward efficient optimization of wind farm layouts: utilizing exact
gradient information.
\newblock {\em Applied Energy}, 179:110--123, 2016.
\newblock \texttt{doi:10.1016/j.apenergy.2016.06.101}.

\bibitem{stanley2019}
A.~P.~J. Stanley, A.~Ning.
\newblock Massive simplification of the wind farm layout optimization problem.
\newblock {\em Wind Energy Science}, 4(4):663--676, 2019.
\newblock \texttt{doi:10.5194/wes-4-663-2019}.

\bibitem{gnegel2021}
F.~Gnegel, A.~F\"ugenschuh, M.~Hagel, S.~Leyffer, M.~Stiemer.
\newblock A solution framework for linear PDE-constrained mixed-integer
problems.
\newblock {\em Mathematical Programming}, 188(2):695--728, 2021.
\newblock \texttt{doi:10.1007/s10107-021-01626-1}.

\bibitem{ott2011fuga}
S.~Ott, J.~Berg, M.~Nielsen.
\newblock Linearised CFD models for wakes.
\newblock Technical Report Ris\o-R-1772(EN), Ris\o\ National Laboratory for
Sustainable Energy, Technical University of Denmark, Roskilde, 2011.

\bibitem{ott2014fuga}
S.~Ott, M.~Nielsen.
\newblock Developments of the offshore wind turbine wake model Fuga.
\newblock Technical Report DTU Wind Energy E-0046, DTU Wind Energy, Roskilde,
2014.

\bibitem{nygaard2020turbopark}
N.~G. Nygaard, S.~T. Steen, L.~Poulsen, J.~G. Pedersen.
\newblock Modelling cluster wakes and wind farm blockage.
\newblock {\em Journal of Physics: Conference Series}, 1618(6):062072, 2020.
\newblock \texttt{doi:10.1088/1742-6596/1618/6/062072}.

\bibitem{pedersen2022turbopark}
J.~G. Pedersen, E.~Svensson, L.~Poulsen, N.~G. Nygaard.
\newblock Turbulence Optimized Park model with Gaussian wake profile.
\newblock {\em Journal of Physics: Conference Series}, 2265(2):022063, 2022.
\newblock \texttt{doi:10.1088/1742-6596/2265/2/022063}.

\bibitem{ti2020mlwake}
Z.~Ti, X.~W. Deng, H.~Yang.
\newblock Wake modeling of wind turbines using machine learning.
\newblock {\em Applied Energy}, 257:114025, 2020.
\newblock \texttt{doi:10.1016/j.apenergy.2019.114025}.

\bibitem{li2023graphwake}
S.~Li, M.~Zhang, M.~D. Piggott.
\newblock End-to-end wind turbine wake modelling with deep graph representation
learning.
\newblock {\em Applied Energy}, 339:120928, 2023.
\newblock \texttt{doi:10.1016/j.apenergy.2023.120928}.

\bibitem{li2022sfnet}
R.~Li, J.~Zhang, X.~Zhao.
\newblock Multi-fidelity modeling of wind farm wakes based on a novel
super-fidelity network.
\newblock {\em Energy Conversion and Management}, 270:116185, 2022.
\newblock \texttt{doi:10.1016/j.enconman.2022.116185}.

\bibitem{bempedelis2024bo}
N.~Bempedelis, F.~Gori, A.~Wynn, S.~Laizet, L.~Magri.
\newblock Data-driven optimisation of wind farm layout and wake steering with
large-eddy simulations.
\newblock {\em Wind Energy Science}, 9(4):869--882, 2024.
\newblock \texttt{doi:10.5194/wes-9-869-2024}.

\bibitem{park2019pgnn}
J.~Park, J.~Park.
\newblock Physics-induced graph neural network: an application to wind-farm
power estimation.
\newblock {\em Energy}, 187:115883, 2019.
\newblock \texttt{doi:10.1016/j.energy.2019.115883}.

\bibitem{li2025hybrid}
P.~Li, Y.~Che, A.~Hua, L.~Wang, M.~Zheng, X.~Guo.
\newblock A data-physics hybrid-driven layout optimization framework for
large-scale wind farms.
\newblock {\em Applied Energy}, 392:125908, 2025.
\newblock \texttt{doi:10.1016/j.apenergy.2025.125908}.

\bibitem{wang2024kriging}
Z.~Wang, Y.~Tu, K.~Zhang, Z.~Han, Y.~Cao, D.~Zhou.
\newblock An optimization framework for wind farm layout design using
CFD-based Kriging model.
\newblock {\em Ocean Engineering}, 293:116644, 2024.
\newblock \texttt{doi:10.1016/j.oceaneng.2023.116644}.

\bibitem{shao2025pce}
Y.-X. Shao, Z.-F. Wang, S.~W. Naung, K.~Zhang, Y.~Yao, D.~Zhou.
\newblock Towards high-fidelity wind farm layout optimization using
polynomial chaos expansion and {Kriging} model.
\newblock {\em Energy}, 338:138820, 2025.
\newblock \texttt{doi:10.1016/j.energy.2025.138820}.

\bibitem{ali2024diffusion}
K.~Ali, T.~Stallard, P.~Ouro.
\newblock A diffusion-based wind turbine wake model.
\newblock {\em Journal of Fluid Mechanics}, 1001:A13, 2024.
\newblock \texttt{doi:10.1017/jfm.2024.1077}.

\bibitem{rethore2014topfarm}
P.-E. R\'ethor\'e, P.~Fuglsang, G.~C. Larsen, T.~Buhl, T.~J. Larsen,
H.~A. Madsen.
\newblock TOPFARM: multi-fidelity optimization of wind farms.
\newblock {\em Wind Energy}, 17(12):1797--1816, 2014.
\newblock \texttt{doi:10.1002/we.1667}.

\bibitem{kajakamaludeen2016}
S.~M. Kaja Kamaludeen, A.~H. van Zuijlen, H.~Bijl.
\newblock Surrogate based wind farm layout optimization using manifold
mapping.
\newblock {\em Journal of Physics: Conference Series}, 753:092005, 2016.
\newblock \texttt{doi:10.1088/1742-6596/753/9/092005}.

\bibitem{song2023surrogate}
D.~Song, X.~Shen, Y.~Gao, L.~Wang, X.~Du, Z.~Xu, Z.~Zhang, C.~Huang,
J.~Yang, M.~Dong, Y.~H. Joo.
\newblock Application of surrogate-assisted global optimization algorithm
with dimension-reduction in power optimization of floating offshore wind
farm.
\newblock {\em Applied Energy}, 351:121891, 2023.
\newblock \texttt{doi:10.1016/j.apenergy.2023.121891}.

\bibitem{kim2026gnn}
T.~Kim.
\newblock A physics-aware surrogate model based on a graph neural network for
mitigating directional discretization bias in wind farm layout optimization.
\newblock Preprint, SSRN, 2026. \texttt{doi:10.2139/ssrn.6923116}.

\bibitem{beale1970}
E.~M.~L.~Beale, J.~A.~Tomlin.
\newblock Special facilities in a general mathematical programming system for
non-convex problems using ordered sets of variables.
\newblock In J.~Lawrence, editor, {\em Proceedings of the Fifth International
Conference on Operational Research}, pages 447--454. Tavistock Publications,
London, 1970.

\bibitem{hexaly}
Hexaly.
\newblock {\em Hexaly Optimizer}, version 15.0 [software], 2026.
\newblock \texttt{https://www.hexaly.com/} (accessed September 2026).

\bibitem{openfoam}
H.~G. Weller, G.~Tabor, H.~Jasak, C.~Fureby.
\newblock A tensorial approach to computational continuum mechanics using
object-oriented techniques.
\newblock {\em Computers in Physics}, 12(6):620--631, 1998.
\newblock \texttt{doi:10.1063/1.168744}.
\newblock Version v2412 (December 2024) used here,
\texttt{https://www.openfoam.com} (accessed September 2026).

\bibitem{vanderlaan2015kefp}
M.~P. van~der Laan, N.~N. S{\o}rensen, P.-E. R\'ethor\'e, J.~Mann, M.~C.
Kelly, N.~Troldborg, et al.
\newblock An improved $k$-$\varepsilon$ model applied to a wind turbine wake
in atmospheric turbulence.
\newblock {\em Wind Energy}, 18(5):889--907, 2015.
\newblock \texttt{doi:10.1002/we.1736}.

\bibitem{vanderlaan2015double}
M.~P. van~der Laan, N.~N. S{\o}rensen, P.-E. R\'ethor\'e, J.~Mann, M.~C.
Kelly, N.~Troldborg.
\newblock The $k$-$\varepsilon$-$f_P$ model applied to double wind turbine
wakes using different actuator disk force methods.
\newblock {\em Wind Energy}, 18(12):2223--2240, 2015.
\newblock \texttt{doi:10.1002/we.1816}.

\bibitem{vanderlaan2015farms}
M.~P. van~der Laan, N.~N. S{\o}rensen, P.-E. R\'ethor\'e, J.~Mann, M.~C.
Kelly, N.~Troldborg, et al.
\newblock The $k$-$\varepsilon$-$f_P$ model applied to wind farms.
\newblock {\em Wind Energy}, 18(12):2065--2084, 2015.
\newblock \texttt{doi:10.1002/we.1804}.

\bibitem{baungaard2022}
M.~Baungaard, S.~Wallin, M.~P. van~der Laan, M.~Kelly.
\newblock Wind turbine wake simulation with explicit algebraic Reynolds
stress modeling.
\newblock {\em Wind Energy Science}, 7(5):1975--2002, 2022.
\newblock \texttt{doi:10.5194/wes-7-1975-2022}.

\bibitem{actuatordiskfoam}
M.~Baungaard.
\newblock {\em actuatorDiskFoam}: a simple actuator disk model for
OpenFOAM [software, MIT licence].
\newblock \texttt{https://github.com/mchba/actuatorDiskFoam} (accessed
November 2026).

\bibitem{gaumond2014}
M.~Gaumond, P.-E. R\'ethor\'e, S.~Ott, A.~Pe\~na, A.~Bechmann, K.~S. Hansen.
\newblock Evaluation of the wind direction uncertainty and its impact on wake
modeling at the Horns Rev offshore wind farm.
\newblock {\em Wind Energy}, 17(8):1169--1178, 2014.
\newblock \texttt{doi:10.1002/we.1625}.

\bibitem{gurobi}
Gurobi Optimization, LLC.
\newblock {\em Gurobi Optimizer Reference Manual}, version 13.0 [software],
2025.
\newblock \texttt{https://www.gurobi.com/} (accessed September 2026).

\bibitem{sood2022lillgrund}
I.~Sood, E.~Simon, A.~Vitsas, B.~Blockmans, G.~C. Larsen, J.~Meyers.
\newblock Comparison of large eddy simulations against measurements from the
Lillgrund offshore wind farm.
\newblock {\em Wind Energy Science}, 7(6):2469--2489, 2022.
\newblock \texttt{doi:10.5194/wes-7-2469-2022}.

\bibitem{glover1997tabu}
F.~Glover, M.~Laguna.
\newblock {\em Tabu Search}.
\newblock Kluwer Academic Publishers, Boston, 1997.
\newblock \texttt{doi:10.1007/978-1-4615-6089-0}.

\bibitem{thomas2022wec}
J.~J. Thomas, S.~McOmber, A.~Ning.
\newblock Wake expansion continuation: multi-modality reduction in the
wind farm layout optimization problem.
\newblock {\em Wind Energy}, 25(4):678--699, 2022.
\newblock \texttt{doi:10.1002/we.2692}.

\end{thebibliography}
\end{document}